\input ./style/arxiv-general.cfg
\documentclass[aos,MSNbibl,dvips]{arximspdf}
\makeatletter
   \@ifpackageloaded{graphicx}{}{\usepackage{graphicx}}
\makeatother

\doi{10.1214/15-AOS1365}
\volume{43}
\issue{6}
\pubyear{2015}
\firstpage{2766}
\lastpage{2794}
\docsubty{FLA}

\makeatletter
\newcommand{\rrvert}{\vert}
\newcommand{\rrVert}{\Vert}
\newcommand{\llvert}{\vert}
\newcommand{\llVert}{\Vert}
\newtheorem{tm}{Theorem}
\newtheorem{la}{Lemma}
\newtheorem{cy}{Corollary}
\newproclaim{co}{Condition}
\newproclaim{rem}{Remark}
\newproclaim{expp}{Experiment}
\makeatother

\begin{document}
\begin{frontmatter}

\title{Estimating the smoothness of a Gaussian random field from
irregularly spaced data via higher-order quadratic variations}
\runtitle{Quadratic variation}

\begin{aug}
\author[A]{\fnms{Wei-Liem}~\snm{Loh}\corref{}\ead[label=e1]{stalohwl@nus.edu.sg}}
\runauthor{W.-L. Loh}
\affiliation{National University of Singapore}
\dedicated{Dedicated to Charles M. Stein on his ninety-fifth birthday}
\address[A]{Department of Statistics and Applied Probability\\
National University of Singapore\\
Singapore 117546\\
Republic of Singapore\\
\printead{e1}}
\end{aug}

%
\received{\smonth{3} \syear{2015}}
%
\revised{\smonth{7} \syear{2015}}

%
\begin{abstract}
This article introduces a method for estimating the smoothness of a stationary,
isotropic Gaussian random field from irregularly spaced data. This
involves novel
constructions of higher-order quadratic variations and the
establishment of the
corresponding fixed-domain asymptotic theory.
In particular, we consider:

\hspace*{5pt}(i) higher-order quadratic variations using nonequispaced line
transect data,

\hspace*{2pt}(ii) second-order quadratic variations from a sample of Gaussian random
field observations taken along a smooth curve in ${\mathbb R}^2$,

(iii) second-order quadratic variations based on deformed lattice data
on~${\mathbb R}^2$.

Smoothness estimators are proposed that are strongly consistent under
mild assumptions.
Simulations indicate that these estimators perform well for moderate
sample sizes.
\end{abstract}

\begin{keyword}[class=AMS]
\kwd[Primary ]{62M30}
\kwd[; secondary ]{62M40}
\end{keyword}
\begin{keyword}
\kwd{Gaussian random field}
\kwd{higher-order quadratic variation}
\kwd{irregularly spaced data}
\kwd{Mat\'{e}rn covariance}
\kwd{smoothness estimation}
\end{keyword}
\end{frontmatter}

\section{Introduction}\label{sec1}

In spatial statistics, it is common practice to model the response as a
realization of a Gaussian random field; cf. \cite{cd1,c1,s1}
and the references cited therein. Let $X$ be a stationary, isotropic
Gaussian random field on ${\mathbb R}^d$ with mean $\mu= E\{ X( {\mathbf
x})\}$ and covariance function
\begin{eqnarray}
\label{eq:1.01} \qquad K( {\mathbf x}, {\mathbf y} ) &=& \operatorname{Cov} \bigl\{
X( {\mathbf x} ), X( {\mathbf y }) \bigr\}
\nonumber
\nonumber
\\[-8pt]
\\[-8pt]
\nonumber
&=& \sum_{j =0}^{\lfloor\nu\rfloor}
\beta_j \llVert{\mathbf x} - {\mathbf y} \rrVert^{2 j} +
\beta_\nu^* G_\nu \bigl( \llVert{\mathbf x} - {\mathbf y}
\rrVert \bigr) + r( {\mathbf x}, { \mathbf y}) \qquad\forall{\mathbf x}, {\mathbf
y}\in{\mathbb R}^d,
\end{eqnarray}
where:
\begin{longlist}[(iii)]
\item[(i)] $\nu>0$, $\beta_\nu^*\neq0$, $\beta_0, \ldots, \beta_{\lfloor\nu
\rfloor}$ are constants and $\lfloor\cdot \rfloor$ denotes the greatest
integer function,

\item[(ii)] $r( {\mathbf x}, {\mathbf y} ) = O( \llVert {\mathbf x} -{\mathbf y} \rrVert ^{2
\nu+ \tau} )$ for some constant $\tau>0$ as $ \llVert {\mathbf x}-{\mathbf
y}\rrVert \rightarrow0$, and
$G_\nu: [0,\infty)\rightarrow{\mathbb R}$ such that $G_\nu(0) = 0$
and for all $s > 0$,
%
%
\begin{equation}
G_\nu( s ) = \cases{ s^{2 \nu}, &\quad$\forall\nu\notin
\mathbb Z$, \vspace*{3pt}
\cr
s^{2 \nu} \log( s), &\quad$\forall\nu\in
\mathbb Z$.} \label{eq:1.19}
\end{equation}
Here, $\llVert\cdot \rrVert $ denotes the Euclidean norm in ${\mathbb R}^d$.
This class of covariance functions is very general and it includes, for
example, the Mat\'{e}rn model\vspace*{1pt} (see Section~\ref{sec2.1})
and the exponential family $\exp(-c \llVert {\mathbf x} - {\mathbf y} \rrVert
^{2 \nu})$ for $\nu\in(0,1)$.
Anderes and Stein \cite{as1}, page 721, observe that $\nu$ is a smoothness parameter in
that $X$ is $j$ times mean square differentiable if and only if $j < \nu$.
$\beta_\nu G_\nu$ is known as the principal irregular term of the
covariance function $K$; cf. \cite{m1,s1}.

The\vspace*{1pt} aim of this article is to estimate $\nu$ using $n$ observations
from a single realization of $X$ within a compact domain $\Delta
\subset{\mathbb R}^d$, $d \in\{1,2\}$.
The estimation of $\nu$ has been addressed in the literature under a
number of different conditions by various authors.
In the case of scattered (possibly nonlattice) data,
Anderes and Stein \cite{as1} ignore the unknown function $r$ in (\ref{eq:1.01}) and use
an approximate likelihood method to estimate $\nu$. However, the
accuracy of the estimate is not addressed there.
Im et~al. \cite{isz1} propose a semiparametric method of estimating the
spectral density (and
hence $\nu$) with irregular observations.
These latter two methods have nonnegligible model biases and also
appear to be analytically intractable under fixed-domain asymptotics.
The latter asymptotics imply that as $n \rightarrow\infty$, the $n$
sites get to be increasingly dense in $\Delta$.
Also in the simulations,
Im et~al. \cite{isz1} consider 200 independent
realizations of a Gaussian random field, whereas this article is
concerned with estimating $\nu$ based
on observations from one realization of the underlying Gaussian random field.

In the case of equispaced data on a line transect,
Hall and Wood \cite{hw1} consider a box-counting estimator while
Constantine and Hall \cite{ch1},
Kent and Wood  \cite{kw1} study estimators of $\nu$ based on process increments.
References \cite{ch1,hw1,kw1} assume that $\nu\in(0, 1)$.
Another example of equally spaced data on a line transect
is \cite{il1} where higher-order quadratic variations are used to
construct a consistent estimate for $\nu$ assuming that $\nu\in(D, D+1)$
for some known integer~$D$.

This article proposes a method for estimating $\nu$ using irregularly
spaced data from a possibly differentiable Gaussian random field. The
methodology involves novel constructions of higher-order quadratic
variations and the establishment of the corresponding fixed-domain
asymptotic theory.
The history of quadratic variations started with \cite{l1}. Since then,
this field has grown dramatically; some examples being \cite
{ap1,bcij1,cgpp1,gl1,kg1}.
Most higher-order quadratic variations in the literature have been
based on observations on a regular grid in ${\mathbb R}^d$; cf. \cite
{cw1,il1} and references cited therein.
An exception is \cite{b1} that deals with second-order quadratic
variations using irregularly spaced Gaussian process observations on a
line transect in $\mathbb R$.
However, from definition (3) of \cite{b1}, we observe that Begyn's
second-order quadratic variations are different
in that they depend explicitly on the smoothness of the process, and
hence cannot be evaluated if $\nu$ is unknown.

The remainder of this article is organised as follows.
Section~\ref{sec2} considers the case $d=1$. Let $\varphi: {\mathbb R}
\rightarrow{\mathbb R}$ be a twice continuously differentiable function
satisfying $\varphi(0) =0$, $\varphi(1)=1$ and $\min_{0\leq s\leq1}
\varphi^{(1)} (s) >0$ where $\varphi^{(1)} (s) = d\varphi(s)/ds$.
Define $t_i = \varphi((i-1)/(n-1))$, $i=1,\ldots, n$.
For $\theta\in\{1,2\}$ and $\ell\in{\mathbb Z}^+$,
novel $\ell$th order quadratic variations $V_{\theta, \ell}$
based on the observations $X( t_i)$, $1\leq i\leq n$, are constructed.
Here, $X$ is a
stationary, isotropic Gaussian random field having covariance function
$K$ as in (\ref{eq:1.01}) with $d=1$.
Under fixed-domain asymptotics, Theorem~\ref{tm:2.1} proves the strong
convergence of $V_{\theta, \ell}/E(V_{\theta, \ell})$ under mild conditions.
It is of interest to note that the asymptotic behavior of $V_{\theta,
\ell}$ is critically dependent on whether the smoothness parameter $\nu
$ is greater or less than the order $\ell$
of the quadratic variation. In Section~\ref{sec2.1}, estimators $\hat{\nu}_{a,
\ell}$, $\hat{\nu}_{a,0}$ and $\hat{\nu}_a$ for $\nu$ are proposed.
For $\nu\leq M < \ell\leq10$ where $M$ is a known constant,
Theorem~\ref{tm:2.2} proves that $\hat{\nu}_{a, \ell}$ is strongly
consistent under the assumptions of Theorem~\ref{tm:2.1}(a).
$\hat{\nu}_a$ is a refinement of $\hat{\nu}_{a,0}$ which in turn can be
thought of as a refinement of $\hat{\nu}_{a, \ell}$.
Table~\ref{tab1} summarises the results of a simulation experiment to gauge the
accuracy of $\hat{\nu}_{a,0}$ and $\hat{\nu}_a$.

For $\theta, \ell\in\{1,2\}$, Section~\ref{sec3} describes the construction of
$\ell$th-order quadratic variations $\tilde{V}_{\theta, \ell}$
from irregularly spaced data taken along a fixed smooth curve $\gamma$
in ${\mathbb R}^2$.
Assuming $\nu\in(0, \ell)$, Theorem~\ref{tm:4.1} proves the strong
convergence of $\tilde{V}_{\theta, \ell}/E( \tilde{V}_{\theta, \ell})$
under weak conditions.
In Section~\ref{sec3.1}, the estimators $\hat{\nu}_{b, 2}$
and $\hat{\nu}_b$ for
$\nu\in(0, 2)$ are proposed.
Theorem~\ref{tm:9.10} shows that $\hat{\nu}_{b,2}$ and $\hat{\nu}_b$
are strongly consistent estimators for $\nu\in(0, 2)$ under the
assumptions of Theorem~\ref{tm:4.1}.
Table~\ref{tab2} reports the accuracy of $\hat{\nu}_{b,2}$ and $\hat{\nu}_b$ in
a simulation experiment where the data are taken along an arc of the
unit circle.

In Section~\ref{sec4}, the Gaussian random field on ${\mathbb R}^2$ is observed
at sites ${\mathbf x}^{i_1, i_2} = \tilde{\varphi} (i_1/n, i_2/n)$, $1\leq
i_1, i_2 \leq n$,
where $\tilde{\varphi}$ is a smooth diffeomorphism.
Second-order quadratic variations $\bar{V}_{\theta, \ell}$, $\theta,
\ell\in\{1, 2\}$, based on the observations\break $X({\mathbf x}^{i_1, i_2})$,
$1\leq i_1, i_2\leq n$, are constructed
where $X$ is a
stationary, isotropic Gaussian random field having covariance function
$K$ as in (\ref{eq:1.01}) with $d=2$.
Theorem~\ref{tm:7.1} establishes conditions where $\bar{V}_{\theta, \ell
}/E ( \bar{V}_{\theta, \ell}) \rightarrow1$ as $n\rightarrow\infty$
almost surely.
Section~\ref{sec4.1} proposes estimators $\hat{\nu}_{c, 1}$ and $\hat{\nu}_{c,
2}$ for $\nu\in(0,2)$.
Theorem~\ref{tm:10.1}
shows that these estimators are strongly consistent under the
assumptions of Theorem~\ref{tm:7.1}. Finally, Table~\ref{tab3} presents the
results of a simulation study on the accuracy of $\hat{\nu}_{c, 1}$ and
$\hat{\nu}_{c,2}$.
Chan and Wood \cite{cw1} consider smoothness estimation of a
nondifferentiable Gaussian random field on ${\mathbb R}^2$ observed on
a regular grid. However, the estimators proposed in \cite{cw1} do not
work for irregularly spaced observations
of a smooth Gaussian random field with $\nu\geq1$.

The \hyperref[appen]{Appendix} and the supplemental article \cite{l2} contain proofs and
related technical results needed in this article.
If $X$ is a locally isotropic, stationary Gaussian random field,
the ideas of this article can still be applied
by choosing the compact domain $\Delta$ suitably small so that $X$ is
close to being isotropic, stationary on $\Delta$.

Likelihood methods are likely to perform well if the data are
relatively sparse and the model is correctly specified.
The situation this article is concerned with is when $r(\cdot,\cdot)$ in (\ref
{eq:1.01}) may not be completely known, the data are relatively dense
and the likelihood is time consuming and difficult (or even impossible)
to compute.

Throughout this article,
$a_n\sim b_n$ denotes $\lim_{n\rightarrow\infty} a_n/b_n = 1$ and
likewise, $a_n\asymp b_n$ denotes $0< \liminf_{n\rightarrow\infty}
a_n/b_n \leq\limsup_{n\rightarrow\infty} a_n/b_n <\infty$.
\end{longlist}

\section{Higher-order quadratic variations using line transect data}\label{sec2}

In this section, the observations of $X$ are taken on a line transect.
Hence, without loss of generality, Section~\ref{sec2} assumes that
$X$ is a Gaussian process having covariance function $K$ as in (\ref
{eq:1.01}) with $d=1$ and
$X(t_{n,1} ), \ldots, X(t_{n,n} )$ are the observed data where
$0 = t_{n, 1} < t_{n, 2} < \cdots< t_{n, n-1} < t_{n, n} = 1$.
For brevity we write, $t_{n,i} = t_i$ and
$X ( t_i) = X_i$, $i =1,\ldots, n$. Define
for $\theta\in\{1,2\}$ and $\ell\in\{ 1, \ldots,
\lfloor(n- 1)/\theta\rfloor\}$
\begin{eqnarray*}
a_{\theta, \ell; i, k} &=& \frac{\ell!}{ \prod_{0 \leq j\leq\ell,
j\neq k} (t_{i+ \theta k} - t_{i + \theta j}) }\qquad\forall k=0,\ldots, \ell,
\nonumber
\\
\nabla_{\theta, \ell} X_i &=& \sum_{k=0}^\ell
a_{\theta, \ell; i, k} X_{i+ \theta k}\qquad\forall i=1,\ldots, n - \theta\ell,
\end{eqnarray*}
and the $\ell$th-order quadratic variation based on $X_1,\ldots, X_n$
to be
%
%
\begin{equation}
V_{\theta, \ell} = \sum_{i=1}^{n- \theta\ell} (
\nabla_{\theta, \ell} X_i )^2. \label{eq:2.1}
\end{equation}

%
\begin{la} \label{la:a.01}
For $\theta\in\{ 1,2 \}$, $\ell\in\{ 1,\ldots, \lfloor(n-
1)/\theta\rfloor\}$ and $i, j \in\{ 1,\ldots, n- \theta\ell\}$, we have
%
%
\begin{eqnarray}
\label{eq:a.1} && \sum_{k=0}^\ell
a_{\theta, \ell; i, k} t_{i+ \theta k}^q = \cases{ 0, &\quad$\forall q=
0, \ldots, \ell-1$, \vspace*{3pt}
\cr
\ell!, &\quad if $q = \ell$,}
\nonumber
\\
&& \sum_{k_1=0}^\ell\sum
_{k_2=0}^\ell a_{\theta, \ell; i, k_1} a_{\theta, \ell; j, k_2 } (
t_{j + \theta k_2} - t_{i+ \theta k_1})^q
\\
&&\qquad = \cases{ 0, &\quad$ \forall q = 0,\ldots, 2 \ell-1$, \vspace*{3pt}
\cr
(-1)^\ell(2 \ell)!, &\quad if $q = 2 \ell$,}
\nonumber
\end{eqnarray}
where we use the convention $0^0 =1$.
\end{la}

The properties of $V_{\theta, \ell}$ rest crucially on the algebraic
identity (\ref{eq:a.1}). Interestingly,
this identity goes way back to \cite{s3}; see also Section~1.2.3,
problem 33, of \cite{k1}. The other result in Lemma~\ref{la:a.01} is a
direct consequence of (\ref{eq:a.1}).
Lemma~\ref{la:a.01} is the reason for using the term ``$\ell$th-order''
in the description of $V_{\theta, \ell}$.
Writing
\begin{eqnarray}
\label{eq:2.4} Y &=& \biggl( \frac{ \nabla_{\theta, \ell} X_1}{\sqrt{ E ( V_{\theta,
\ell} )}}, \ldots, \frac{ \nabla_{\theta, \ell} X_{n- \theta\ell}}{
\sqrt{ E (V_{\theta, \ell} )}}
\biggr)',
\nonumber
\nonumber
\\[-8pt]
\\[-8pt]
\nonumber
\Sigma&=& (\Sigma_{i,j})_{( n- \theta\ell) \times(n - \theta\ell
)} = E \bigl( Y
Y' \bigr),
\end{eqnarray}
we obtain $V_{\theta, \ell}/E ( V_{\theta, \ell} ) = Y' Y = Z' \Sigma
Z$ where $Z \sim N_{n- \theta\ell} (0, I)$.
It follows from \cite{hw2} that
%
%
\begin{eqnarray}
\label{eq:2.5} && P \biggl( \biggl\llvert\frac{ V_{\theta, \ell} }{ E (V_{\theta, \ell}
) } -1 \biggr\rrvert\geq
\varepsilon \biggr)
\nonumber
\\
&&\qquad = P \bigl( \bigl\llvert Z' \Sigma Z - E \bigl(
Z' \Sigma Z \bigr) \bigr\rrvert\geq\varepsilon \bigr)
\\
&&\qquad \leq 2 \exp \biggl\{- C\min \biggl( \frac{\varepsilon}{ \llVert
\Sigma_{\mathrm{abs}}\rrVert _2},
\frac{\varepsilon^2}{ \llVert \Sigma
_{\mathrm{abs}}\rrVert _F^2} \biggr) \biggr\}\qquad\forall\varepsilon>0,
\nonumber
\end{eqnarray}
where $C>0$ is an absolute constant, $\Sigma_{\mathrm{abs}}$ is the $(n-
\theta\ell)\times(n- \theta\ell)$ matrix
with elements $\llvert \Sigma_{i, j}\rrvert $, $i,j=1,\ldots, n -
\theta\ell$, and $\llVert\cdot \rrVert _2$, $\llVert\cdot \rrVert _F$ are the
spectral, Frobenius matrix norms, respectively.

Lemma~\ref{la:a.01} and (\ref{eq:2.5}) are needed in the proof of
Theorem~\ref{tm:2.1} below. The latter provides a way for estimating
$\nu$.
It is convenient to write $\tilde{K}( x-y) = K(x,y)$ for all $x, y\in
\mathbb R$ and $\tilde{K}^{(2 \ell)}$ as the $2 \ell$th derivative of
$\tilde{K}$ if the latter exists.

\begin{co}\label{co1}
For $n\geq2$, define $t_i = \varphi((i-1)/(n-1))$,
$i=1,\ldots, n$, where $\varphi: {\mathbb R} \rightarrow{\mathbb R}$
is a twice continuously differentiable function
satisfying $\varphi(0) =0$, $\varphi(1)=1$ and $\min_{0\leq s\leq1}
\varphi^{(1)} (s) >0$.
\end{co}

We then say that the $t_i$'s are generated by $\varphi$.
It follows from Condition~\ref{co1} that there exist constants $C_{1,0}$ and
$C_{1,1}$ such that
\[
0 < C_{1,0}/n \leq\min_{1 \leq i\leq n-1} ( t_{i+1} -
t_i ) \leq\max_{1\leq i\leq n-1} ( t_{i+1} -
t_i ) \leq C_{1,1}/n. \label{eq:2.6}
\]

Writing $G_\nu(\cdot)$ as in (\ref{eq:1.19}), for $\theta\in\{1,2\}$ and
$\ell\in\{1, \ldots, \lfloor(n-1)/\theta\rfloor\}$, define
%
%
\begin{eqnarray}
\label{eq:2.101} f_{\theta, \ell} ( \nu) = 2 \beta_\nu^* \sum
_{i=1}^{n- \theta\ell} \sum
_{0\leq k_1 < k_2 \leq\ell} a_{\theta, \ell; i, k_1} a_{\theta, \ell; i, k_2} G_\nu(
t_{i+ \theta
k_2} - t_{i+ \theta k_1} )
\nonumber
\\[-8pt]
\\[-8pt]
\eqntext{\forall\nu\in(0, \ell).}
\end{eqnarray}

%
\begin{tm} \label{tm:2.1}
Let $V_{\theta, \ell}$ be as in (\ref{eq:2.1}), $\theta\in\{1,2\}$
and Condition~\ref{co1} holds.

\begin{longlist}[(a)]
\item[(a)] Suppose $\nu\in(0, \ell)$ and $\ell\in\{1, \ldots, \lfloor
(n-1)/\theta\rfloor\wedge10 \}$.
Then
$E( V_{\theta, \ell}) \sim f_{\theta, \ell} ( \nu) \asymp n^{2\ell+1 -
2 \nu}$.
In addition, if $\tilde{K}^{(2 \ell)} (\cdot)$ is a continuous function on
an open interval containing $( 0, 1]$
such that $\llvert \tilde{K}^{(2\ell)} (t)\rrvert \leq C_\ell t^{2 \nu
-2\ell}$, $\forall t\in(0, 1]$, for some constant $C_\ell$, then
\[
\operatorname{Var} \bigl\{ V_{\theta, \ell}/E (V_{\theta, \ell} ) \bigr\} =
\cases{ O \bigl( n^{-1} \bigr), &\quad if $\nu< (4 \ell-1)/4$,
\vspace*{3pt}
\cr
O \bigl\{ n^{-1} \log(n) \bigr\}, &\quad if $\nu= (4
\ell-1)/4$, \vspace*{3pt}
\cr
O \bigl( n^{-4 \ell+ 4 \nu} \bigr), &\quad if $\nu> (4
\ell-1)/4$,}
\]
and $V_{\theta, \ell}/E( V_{\theta, \ell}) \rightarrow1$ almost surely
as $n\rightarrow\infty$.

\item[(b)] Suppose $\nu= \ell$. Then
$E( V_{\theta, \ell}) \sim(-1)^{\ell+1} \beta_\ell^* (2 \ell)! n \log(n)$.
In addition, if $\tilde{K}^{(2 \ell)} (\cdot)$ is a continuous function on
an open interval containing $( 0, 1]$
such that $\llvert \tilde{K}^{(2\ell)} (t)\rrvert \leq C_\ell\log(
2/t )$, $\forall t\in(0, 1]$, for some constant $C_\ell$, then
$\operatorname{Var}\{ V_{\theta, \ell}/E (V_{\theta, \ell}) \} = O\{ \log^{-2} (n) \}$
as $n\rightarrow\infty$.

\item[(c)] Suppose $\nu> \ell$. Then
$E( V_{\theta, \ell}) \sim(-1)^\ell\beta_\ell(2 \ell)! n$.
In addition, if $\tilde{K}^{(2 \ell)} (\cdot)$ is a continuous and not
identically $0$ function on an open interval containing $( 0, 1]$, then
\[
\liminf_{n \rightarrow\infty} \operatorname{Var}\bigl\{ V_{\theta, \ell}/E
(V_{\theta, \ell
}) \bigr\} > 0.
\]
\end{longlist}
\end{tm}

\begin{rem}\label{rem1}
In Theorem~\ref{tm:2.1}(a), we restrict $\ell\leq10$
as we think this will suffice in practical situations.
However, if $V_{\theta, \ell}$ for some $\ell>10$ is required, then as
in the proof of Lemma~\ref{la:a.1}, one need only use, say, \textit{Mathematica}, to verify that $H_\ell(\nu)\neq0$ for all $0\leq\nu
\leq\ell$.
\end{rem}

\subsection{Smoothness estimation on a line transect}\label{sec2.1}

Suppose $0<\nu\leq M < \ell\leq10$ for some known constant $M$.
Under the conditions of Theorem~\ref{tm:2.1}(a), it is easy to
construct a strongly consistent estimator for $\nu$. For example,
taking $\theta=1$,
$\{ 2 \ell+1 - \log(V_{1, \ell})/\log(n)\}/2$
is one such estimator for $\nu$. However, the bias is of order $1/\log
(n)$ which makes it unsuitable for use in practice.
With the notation of Theorem~\ref{tm:2.1}, define $F_{\ell, n}:[0, M]
\rightarrow[0,\infty)$ by
$F_{\ell, n} (\nu^*) = f_{2, \ell} (\nu^*)/f_{1, \ell} (\nu^*)$ for $\nu
^* \in(0, M]$ and $F_{\ell, n} (0) = \lim_{\delta\rightarrow0+}
f_{2, \ell} (\delta)/f_{1, \ell} (\delta)$.
The motivation for taking the ratio on the right-hand side is to
eliminate the nuisance parameter $\beta_{\nu^*}^*$.
We observe from Lemma~\ref{la:a.4}
that $F_{\ell, n} (\cdot)$ is a continuous function.
We shall now construct another strongly consistent estimator $\hat{\nu
}_{a, \ell}$ for $\nu$.
Let $\hat{\nu}_{a, \ell} \in[0, M]$ satisfy
%
%
\begin{equation}
\biggl\{ \frac{ V_{1, \ell} F_{\ell, n} (\hat{\nu}_{a, \ell}) }{ V_{2,
\ell}} - 1 \biggr\}^2 = \min
_{0 \leq\nu^* \leq M} \biggl\{ \frac{ V_{1, \ell} F_{\ell, n} (\nu^*)
}{ V_{2, \ell}} - 1 \biggr\}^2.
\label{eq:2.1.4}
\end{equation}

%
\begin{tm} \label{tm:2.2}
Let $0< \nu\leq M < \ell\leq10$, $\hat{\nu}_{a, \ell}$ be as in (\ref
{eq:2.1.4}) and that
the conditions of Theorem~\ref{tm:2.1}\textup{(a)} are satisfied.
Then $\hat{\nu}_{a, \ell} \rightarrow\nu$ as $n\rightarrow\infty$
almost surely.
\end{tm}

In practice, the upper bound $M$ is usually conservative and can be
significantly larger than the unknown $\nu$.
Thus, $\hat{\nu}_{a, \ell}$ may not perform well for small to moderate
sample size $n$.
Next, we propose an alternative estimator $\hat{\nu}_a$ for $\nu$ that
refines on $\hat{\nu}_{a, \ell}$ by first estimating an interval of
unit width which contains $\nu$.
The algorithm below is motivated by the fact that
$ V_{1, l} F_{l, n} (\nu)/V_{2, l} \rightarrow1$ as $n\rightarrow
\infty$ almost surely if $\nu< l$ (cf. Theorem~\ref{tm:2.1}).

\begin{longlist}[\textit{Step} 2.1.]
\item[\textit{Step} 2.1.] For each $l = 1, \ldots, \lfloor M \rfloor+2$, let
$\tilde{\nu}_{a, l} \in[0, M]$ be such that
\[
\biggl\{ \frac{ V_{1, l} F_{l,n} (\tilde{\nu}_{a, l}) }{ V_{2, l}} - 1 \biggr\}^2 = \min
_{0\leq\nu^* \leq M \wedge l } \biggl\{ \frac{ V_{1, l} F_{l, n} (\nu
^*) }{V_{2, l}} - 1 \biggr\}^2.
\]

\item[\textit{Step} 2.2.] Let $\hat{\nu}_{a, 0}$ be the value of $\tilde{\nu}_{a,
l}$ that minimises $( \tilde{\nu}_{a, l} - \tilde{\nu}_{a, {l +1}}
)^2$, $l = 1, \ldots, \lfloor M \rfloor+1$.
The purpose of $\hat{\nu}_{a, 0}$ is to estimate an interval of unit
width containing $\nu$.

\item[\textit{Step} 2.3.] We improve on $\hat{\nu}_{a, 0}$ by defining the
estimator $\hat{\nu}_a$ for $\nu$ to be
$\hat{\nu}_a = \tilde{\nu}_{a, l}$ where $l$ is the smallest integer
satisfying $l > \hat{\nu}_{a, 0} + 1/4$.
\end{longlist}

\begin{rem}\label{rem2}
The motivation behind step 2.2 is that if $0< \nu< l$,
then $\tilde{\nu}_{a, l} \rightarrow\nu$ and
$( \tilde{\nu}_{a, l} - \tilde{\nu}_{a, {l +1}} )^2 \rightarrow0$ as
$n\rightarrow\infty$ almost surely.
The rationale for step 2.3 is to use $V_{\theta, l}$, $\theta\in\{
1,2\}$ with the smallest integer $l > \nu+ 1/4$ to estimate $\nu$
since it follows from Theorem~\ref{tm:2.1}(a) that
$\operatorname{Var}\{V_{\theta, l}/E(V_{\theta, l})\} = O(n^{-1})$ as $n\rightarrow
\infty$.

As noted by Stein \cite{s1}, a class of covariance functions which has
considerable practical value is
the Mat\'{e}rn class:
%
%
\begin{equation}
K_{\mathrm{Mat}} ( {\mathbf x}, {\mathbf y} ) = \frac{ \sigma^2 (\alpha\llVert
{\mathbf x} - {\mathbf y} \rrVert )^\nu}{ 2^{\nu-1} \Gamma(\nu) } {\mathcal
K}_\nu \bigl(\alpha\llVert{\mathbf x} - {\mathbf y} \rrVert \bigr)
\qquad \forall{\mathbf x}, {\mathbf y} \in{\mathbb R}^d,
\label{eq:2.1.1}
\end{equation}
where $\nu, \alpha, \sigma$ are strictly positive constants and
${\mathcal K}_\nu$ is the modified
Bessel function of the second kind.
This implies that if $\nu$ is not an integer,
%
%
\begin{eqnarray}
\label{eq:2.1.2} K_{\mathrm{Mat}} ( {\mathbf x}, {\mathbf y}) &=&
\sigma^2 \sum_{k=0}^\infty
\frac{ \alpha^{2 k} \llVert {\mathbf x} - {\mathbf y} \rrVert
^{2 k} }{ 2^{2k} k! \prod_{i=1}^k (i - \nu) }
\nonumber
\\[-8pt]
\\[-8pt]
\nonumber
&&{} - \frac{ \pi\sigma^2}{
\Gamma(\nu) \sin(\nu\pi) } \sum_{k=0}^\infty
\frac{ \alpha^{2k+ 2 \nu} \llVert {\mathbf x} - {\mathbf y}
\rrVert ^{2 k+ 2 \nu} }{ 2^{2k+ 2 \nu} k! \Gamma(k+1 +\nu)}.
\end{eqnarray}
On the other hand, if $\nu\in\{ 1,2,\ldots\}$,
%
%
\begin{eqnarray}
\label{eq:2.1.3} K_{\mathrm{Mat}} ( {\mathbf x}, {\mathbf y}) &=&
\frac{ 2 \sigma^2 }{ (\nu-1)! } \Biggl\{ (-1)^{\nu+1} \log \biggl( \frac{ \alpha\llVert {\mathbf x} - {\mathbf y}
\rrVert }{ 2}
\biggr) \sum_{k=0}^\infty\frac{ ( \alpha\llVert {\mathbf x} - {\mathbf y} \rrVert /2)^{
2 (\nu+ k) } }{ k! (\nu+ k )! }
\nonumber
\\
&&{} + \frac{1}{2} \sum_{k=0}^{\nu-1}
(-1)^k \frac{ (\nu-k-1 )!}{k!} \biggl(\frac{ \alpha\llVert {\mathbf x} -
{\mathbf y} \rrVert }{2}
\biggr)^{2 k}
\\
&&{} + \frac{(-1)^\nu}{2} \sum_{k=0}^\infty
\bigl[ \psi(k+1) + \psi( \nu+k +1) \bigr] \frac{ ( \alpha\llVert
{\mathbf x} - {\mathbf y} \rrVert /2)^{2 (\nu+ k )} }{k! (\nu+ k)!} \Biggr\},
\nonumber
\end{eqnarray}
where $\psi(\cdot)$ is the digamma function, that is,
$\psi(k) = -\gamma+ \sum_{i=1}^{k-1} i^{-1}$, $k= 1,2, \ldots,$ and
$\gamma= 0.5772 \ldots$ is Euler's constant.
Equations (\ref{eq:2.1.2}) and (\ref{eq:2.1.3}) indicate that $K_{\mathrm{Mat}} (\cdot,\cdot)$
can be expressed as in (\ref{eq:1.01}).
It can be easily shown that $\tilde{K}_{\mathrm{Mat}}^{(2 \ell)}(\cdot)$
satisfies the conditions in Theorem~\ref{tm:2.1} as well.

In order to gauge the finite sample accuracy of $\hat{\nu}_{a, 0}$ and
$\hat{\nu}_a$, a simulation experiment is conducted.
In this experiment, $X$ is a stationary Gaussian process having mean 0
and Mat\'{e}rn covariance function $K_{\mathrm{Mat}}$ as in (\ref{eq:2.1.1})
with $\sigma= \alpha= d = 1$.
All computations are performed using the software {\em Mathematica}.
\end{rem}

\begin{expp}\label{expp1}
Set $n = 200$, $M = 2.5$ and $\nu= 0.1$, 0.5, 0.9,
1, 1.1, 1.5, 1.9, 2, 2.1, 2.5.
For $i=1,\ldots, n$, $X_i = X (t_i)$ where $t_i = \varphi
((i-1)/(n-1))$ with $\varphi( s) = s(s+1)/2$ for all $s \in\mathbb R$.
The mean absolute errors of $\hat{\nu}_{a, 0}$ and $\hat{\nu}_a$ are
computed over 100 replications.

%
\begin{table}[b]
\tabcolsep=0pt
\tablewidth=215pt
\caption{Experiment~\protect\ref{expp1}. Simulation results in estimating $\nu$
using data generated by $\varphi(s) = s(s+1)/2$ over 100 replications
(standard errors within parentheses)}\label{tab1}
\begin{tabular*}{\tablewidth}{@{\extracolsep{\fill}}@{}lcc@{}}
\hline
$\bolds{\nu}$ & $\bolds{E\llvert \hat{\nu}_{a, 0} - \nu\rrvert }$ & $\bolds{E\llvert \hat {\nu}_a - \nu\rrvert }$ \\
\hline
$0.1$ 
& 0.109 (0.010) & 0.061 (0.004) \\
$0.5$ 
& 0.128 (0.013) & 0.057 (0.005) \\
$0.9$ 
& 0.107 (0.011) & 0.074 (0.006) \\
$1.0$ 
& 0.108 (0.009) & 0.074 (0.006) \\
$1.1$ 
& 0.113 (0.008) & 0.069 (0.005) \\
$1.5$ 
& 0.093 (0.008) & 0.052 (0.005) \\
$1.9$ 
& 0.077 (0.005) & 0.078 (0.005) \\
$2.0$ 
& 0.082 (0.005) & 0.069 (0.005) \\
$2.1$ 
& 0.079 (0.006) & 0.063 (0.006) \\
$2.5$ 
& 0.053 (0.006) & 0.049 (0.004) \\ \hline
\end{tabular*}
\end{table}

The results from Experiment~\ref{expp1} are presented in Table~\ref{tab1}. The mean
absolute errors of $\hat{\nu}_{a,0}$ and
$\hat{\nu}_a$ indicate that $\hat{\nu}_a$ is significantly more
accurate than $\hat{\nu}_{a, 0}$ thus vindicating step 2.3.

It is well known that simulating a stationary Gaussian process on
$[0,1]$ when $n$ and $\nu$ are large is a difficult problem; cf. \cite{s2,wc1}.
This is especially so if the data are irregularly spaced. This is the
reason why we set the upper limit for the true value of $\nu$ to be 2.5
in Experiment~\ref{expp1}.
For larger values of $\nu$, the {\em Mathematica} routine
MultinormalDistribution$[\cdot,\cdot]$ returns with the error message that the
covariance matrix is not sufficiently positive definite to complete the
Cholesky decomposition to reasonable accuracy.
Finally, we have not compared the performance of $\hat{\nu}_a$ to the
estimator proposed in \cite{as1} because the latter estimator requires a
choice of neighbourhood blocks and there are no formal guidelines given
for choosing these blocks.
\end{expp}

\section{Second-order quadratic variation along a curve in ${\mathbb R}^2$ using 3 or more points}\label{sec3}

In this section, the observations of $X$ are taken along a fixed curve
$\gamma$ in ${\mathbb R}^2$ where
$X$ is a Gaussian random field having covariance function $K$ as in
(\ref{eq:1.01}) with $d = 2$.
More precisely, we assume
the following.

\begin{co}\label{co2}
There exist strictly positive constants $\varepsilon
, L$ such that $\gamma: (-\varepsilon, L+\varepsilon) \rightarrow
{\mathbb R}^2$ is
a $C^2$-curve parametrised by its arc length.
In\vspace*{1pt} particular
writing $\gamma(t) = (\gamma_1(t), \gamma_2(t))'$
and its $k$th derivative by $\gamma^{(k)} (t) = (\gamma_1^{(k)} (t),\break
\gamma_2^{(k)} (t))'$, we have (i) $\gamma^{(2)} (t)$ exists and is continuous,
and (ii) $\llVert \gamma^{(1)} (t)\rrVert = 1$ for all $t\in
(-\varepsilon, L+\varepsilon)$.
\end{co}

\begin{co}\label{co3}
There exists a constant $C_3 >0$ such that
\[
\bigl\llVert\gamma \bigl(t^* \bigr) - \gamma(t) \bigr\rrVert\geq
C_3 \bigl\llvert t^* - t \bigr\rrvert\qquad\forall t^*, t \in [0,
L].
\]
\end{co}

\begin{co}\label{co4}
For $n\geq2$, $t_{n,i} = \varphi( L (i-1)/(n-1))$,
$i=1,\ldots, n$, where $\varphi: {\mathbb R} \rightarrow{\mathbb R}$
is a twice continuously differentiable function satisfying $\varphi(0)
=0$, $\varphi(L)=L$ and $\min_{0\leq s\leq L} \varphi^{(1)} (s) >0$.
\end{co}

\begin{rem}\label{rem3}
Condition~\ref{co3} ensures that the curve $\gamma$ is
reasonably well behaved, for example, $\gamma$ does not intersect itself.

In this section, we shall write $t_i = t_{n,i}$, $X ( \gamma(t_i )) =
X_i$ and $d_{i,j} = \llVert \gamma(t_i) - \gamma(t_j ) \rrVert $
for all $1 \leq i, j \leq n$.
$X_1,\ldots, X_n$ represent the observations of $X$ that are made on
the curve $\gamma$.
Define for $\theta, \ell\in\{1,2\}$,
\begin{eqnarray*}
b_{\theta, \ell; i, k} &=& \frac{ \ell}{ \prod_{ 0 \leq j\leq\ell, j
\neq k } ( d_{i, i+ \theta k} - d_{i, i+ \theta j} ) }\qquad\forall k= 0, \ldots, \ell,
\nonumber
\\
\tilde{\nabla}_{\theta, \ell} X_i &=& \sum
_{k= 0}^\ell b_{\theta, \ell; i, k} X_{i+ \theta k}
\qquad\forall i =1,\ldots, n- \theta\ell,
\end{eqnarray*}
and the $\ell$th-order quadratic variation based on $X_1,\ldots, X_n$
to be
%
%
\begin{equation}
\tilde{V}_{\theta, \ell} = \sum_{i = 1}^{n - \theta\ell}
( \tilde{\nabla}_{\theta, \ell} X_i )^2.
\label{eq:4.1}
\end{equation}
Then we observe from
Lemma~\ref{la:a.01} that for $\theta, \ell\in\{1,2\}$,
%
%
\begin{eqnarray}
\sum_{k=0}^\ell b_{\theta, \ell; i, k}
d^q_{i, i+\theta k} &=& \cases{ 0, &\quad if $q=0, \ldots, \ell-1$,
\vspace*{3pt}
\cr
\ell, &\quad if $q= \ell$.} \label{eq:4.2}
\end{eqnarray}
Equation (\ref{eq:4.2}) is the motivation for calling $\tilde{V}_{\theta, \ell}$
``$\ell$th-order''.
Define $\tilde{K} ({\mathbf x} - {\mathbf y} ) = K ({\mathbf x}, {\mathbf y})$ for all
${\mathbf x}, {\mathbf y} \in{\mathbb R}^2$. Letting $a_1$ and $a_2$ be
nonnegative integers, we write ${\mathbf x} = (x_1, x_2)'$ and
\[
\tilde{K}^{(a_1, a_2)} ({\mathbf x}) = \frac{\partial^{a_1+ a_2}}{\partial
x_1^{a_1} \,\partial x_2^{a_2}} \tilde{K} ( {\mathbf
x}).
\]
Let $\Gamma$ be a compact, convex set in ${\mathbb R}^2$ such that $\{
\gamma(t): 0\leq t\leq L\} \subset\Gamma$.
For $\theta, \ell\in\{1,2\}$, define
%
%
\begin{eqnarray}
\label{eq:2.100} \tilde{f}_{\theta, \ell} ( \nu) = 2 \beta_\nu^* \sum
_{i=1}^{n- \theta\ell} \sum
_{0\leq k_1 < k_2 \leq\ell} b_{\theta, \ell; i, k_1} b_{\theta, \ell; i, k_2} G_\nu(
d_{i+ \theta
k_1, i+ \theta k_2} )
\nonumber
\\[-8pt]
\\[-8pt]
\eqntext{\forall\nu\in(0, \ell). }
\end{eqnarray}
\end{rem}

%
\begin{tm} \label{tm:4.1}
Let $\theta, \ell\in\{1,2\}$, $\nu\in(0, \ell)$ and $\tilde
{V}_{\theta, \ell}$ be as in (\ref{eq:4.1}). Suppose Conditions~\ref{co2},~\ref{co3}
and~\ref{co4} hold. Then
$E( \tilde{V}_{\theta, \ell}) \sim\tilde{f}_{\theta, \ell} (\nu)
\asymp n^{ 2 \ell+1 - 2 \nu}$.
In addition, suppose \textup{(i)} $\tilde{K}^{(a_1, a_2)} (\cdot)$ is a continuous
function on an open set containing $\{ {\mathbf x}-{\mathbf y}: {\mathbf x}, {\mathbf
y}\in\Gamma, {\mathbf x}\neq{\mathbf y}\}$
whenever $a_1, a_2$ are nonnegative integers such that $a_1 + a_2= 2
\ell$ and
\textup{(ii)}~there exists a constant $C_4$ such that
\[
\bigl\llvert\tilde{K}^{(a_1, a_2)} ( {\mathbf x} - {\mathbf y}) \bigr\rrvert
\leq C_4 \llVert{\mathbf x} - {\mathbf y} \rrVert^{2 \nu- 2 \ell}
\qquad \forall2 \leq a_1+ a_2 \leq2 \ell,
\]
for all ${\mathbf x}, {\mathbf y} \in\Gamma$, ${\mathbf x}\neq{\mathbf y}$.
Then $\tilde{V}_{\theta, \ell}/E( \tilde{V}_{\theta, \ell}) \rightarrow
1$ almost surely
and
\[
\operatorname{Var} \bigl\{ \tilde{V}_{\theta, \ell}/E ( \tilde{V}_{\theta, \ell})
\bigr\} = \cases{ O \bigl( n^{-1} \bigr), &\quad if $\nu< (4
\ell-1)/4$, \vspace*{3pt}
\cr
O \bigl\{ n^{-1} \log(n) \bigr\}, & \quad
if $\nu= (4 \ell-1)/4$, \vspace*{3pt}
\cr
O \bigl( n^{- 4 \ell+ 4 \nu} \bigr), &\quad
if $\nu> (4\ell-1)/4$,}\qquad\mbox {as $n \rightarrow\infty$}.
\]
\end{tm}

The proof of Theorem~\ref{tm:4.1} is can be found in the
supplemental article \cite{l2}.

\subsection{Estimating $\nu$ along a smooth curve in ${\mathbb R}^2$}\label{sec3.1}

Let $\theta, \ell\in\{1,2\}$ and $\tilde{f}_{\theta, \ell}$ be as in
(\ref{eq:2.100}).
Define $\tilde{F}_{\ell, n}: [0, \ell] \rightarrow{\mathbb R}$ by
$\tilde{F}_{\ell, n} (\nu^*) =
\tilde{f}_{2, \ell} (\nu^*)/ \tilde{f}_{1, \ell} (\nu^*)$ $\forall\nu
^* \in(0, \ell)$, $\tilde{F}_{\ell, n} (0) = \lim_{\delta\rightarrow
0+} \tilde{f}_{2, \ell} ( \delta)/\tilde{f}_{1, \ell} (\delta)$ and
$\tilde{F}_{\ell, n} ( \ell) = \lim_{\delta\rightarrow0+} \tilde
{f}_{2, \ell} (\ell- \delta)/\tilde{f}_{1, \ell} ( \ell- \delta)$. We
observe from Lemma~\ref{la:a.4} that $\tilde{F}_{\ell, n} (\cdot)$ is a
continuous function.
Let $\hat{\nu}_{b, \ell} \in[0, \ell]$ satisfy
%
%
\begin{equation}
\biggl\{ \frac{ \tilde{V}_{1, \ell} \tilde{F}_{\ell, n} (\hat{\nu}_{b,
\ell}) }{ \tilde{V}_{2, \ell}} - 1 \biggr\}^2 = \min
_{0 \leq\nu^* \leq\ell} \biggl\{ \frac{ \tilde{V}_{1, \ell} \tilde
{F}_{\ell, n} (\nu^*) }{ \tilde{V}_{2, \ell}} - 1 \biggr\}^2.
\label{eq:9.69}
\end{equation}

%
\begin{tm} \label{tm:9.10}
Let $\ell\in\{1,2\}$, $\nu\in(0,\ell)$ and $\hat{\nu}_{b, \ell}$ be
as in (\ref{eq:9.69}). Suppose
the conditions of Theorem~\ref{tm:4.1} are satisfied.
Then $\hat{\nu}_{b, \ell} \rightarrow\nu$ as $n\rightarrow\infty$
almost surely.
\end{tm}

Using Lemma~\ref{la:3.9}, the proof of Theorem~\ref{tm:9.10} is similar
to that of Theorem~\ref{tm:2.2} and will be omitted.
Theorem~\ref{tm:9.10} proves that $\hat{\nu}_{b,2}$ is a strongly
consistent estimator for $\nu\in(0, 2)$.
However, if $\nu$ is close to 0, the upper bound of 2 is conservative
and a better estimator for $\nu$ is $\hat{\nu}_{b,1}$.
Consequently, we propose the following alternative estimator $\hat{\nu
}_b$ for $\nu\in(0,2)$.

Compute $\hat{\nu}_{b,2}$. If $\hat{\nu}_{b,2} > 3/4$, define $\hat{\nu
}_b = \hat{\nu}_{b,2}$. If $\hat{\nu}_{b,2} \leq3/4$,
compute $\hat{\nu}_{b,1}$ and define $\hat{\nu}_b = \hat{\nu}_{b,1}$.
$3/4$ is motivated by the last statement of Theorem~\ref{tm:4.1} when
\mbox{$\ell=1$}. Under the assumptions of Theorem~\ref{tm:9.10}, we observe
that $\hat{\nu}_b$ is also a strongly consistent
estimator for $\nu$.

\begin{rem}\label{rem4}
We observe from the definitions of $\hat{\nu}_{b,2},
\hat{\nu}_b$ that $\varphi$ need not be explicitly known; only the
bijection $\Phi: \{1,\ldots, n\} \rightarrow\{ \phi_i = \varphi
(L(i-1)/(n-1)): i=1,\ldots, n\}$ is required.
There are 2 feasible bijections, namely $i\mapsto\varphi(
L(i-1)/(n-1))$ and $i\mapsto\varphi( L(n-i)/(n-1))$. Either bijection
will do.
If the curve $\gamma$ is known, then $\Phi$ can be recovered from the
ordering of the $\phi_i$'s along $\gamma$. On the other hand if $\gamma
$ is unknown, then $\Phi$ can be obtained
using the fact that ``adjacent'' $\phi_i$'s are closest to each other
for sufficiently large $n$, that is, $\max_{k=-1, 1} \llVert \phi
_{j+k} -\phi_j\rrVert < \min_{i\neq j-1, j, j+1} \llVert \phi_i - \phi
_j\rrVert $.
This greatly increases the utility of $\hat{\nu}_{b,2}, \hat{\nu}_b$ in
applications. The following is an algorithm for determining $\Phi$ when
$\gamma$ is unknown:

\begin{longlist}[\textit{Step} 3.1.]
\item[\textit{Step} 3.1.] First, arbitrarily pick an element $y_0$ from the set
$\{\phi_1,\ldots, \phi_n\}$.

\item[\textit{Step} 3.2.] Let $y_1$ denote the element of $\tilde{Y}_0 = \{\phi
_1,\ldots, \phi_n\}\setminus\{y_0\}$ such that $\llVert
y_1-y_0\rrVert = \min\{ \llVert \phi_i-y_0\rrVert: \phi_i\in\tilde
{Y}_0\}$.

\item[\textit{Step} 3.3.] Given $Y_{-k, l} = \{ y_{-k}, \ldots, y_0, \ldots, y_l\}
$, let $y^*, y^{**}$ denote elements of $\tilde{Y}_{-k,l} = \{\phi
_1,\ldots, \phi_n\}\setminus Y_{-k,l}$ such that $\llVert y^*
-y_{-k}\rrVert = \min\{ \llVert \phi_i-y_{-k}\rrVert: \phi_i\in\tilde
{Y}_{-k,l}\}$ and
$\llVert y^{**}-y_l\rrVert = \min\{ \llVert \phi_i-y_l\rrVert: \phi
_i\in\tilde{Y}_{-k,l}\}$. If $\llVert y^* -y_{-k}\rrVert \leq
\llVert y^{**}-y_l\rrVert $, define $y_{-(k+1)} = y^*$ and if
$\llVert y^* -y_{-k}\rrVert > \llVert y^{**}-y_l\rrVert $, define
$y_{l+1} = y^{**}$.

\item[\textit{Step} 3.4.] Repeat step 3.3 until the cardinality of $Y_{-k,l}$ is
$n$. The ordering of this $Y_{-k,l}$, with cardinality $n$ and
$l=n-k-1$, gives the required bijection $\Phi(i) = y_{-k+i-1}$,
$i=1,\ldots, n$.
\end{longlist}

A simulation experiment is conducted to gauge the finite sample
accuracy of $\hat{\nu}_{b,2}$ and $\hat{\nu}_b$.
In Experiment~\ref{expp2} below:
\begin{longlist}[(ii)]
\item[(i)] $X$ is a stationary Gaussian random field having mean 0 and Mat\'
{e}rn covariance function $K_{\mathrm{Mat}}$ as in (\ref{eq:2.1.1}) with
$\sigma= \alpha= 1$ and $d = 2$.

%
\begin{table}[t]
\tabcolsep=0pt
\tablewidth=212pt
\caption{Experiment \protect\ref{expp2}. Simulation results in estimating $\nu$
with nonequispaced data along a curve $\gamma$ over 100 replications
(standard error within parentheses)}\label{tab2}
\begin{tabular*}{\tablewidth}{@{\extracolsep{\fill}}@{}lcc@{}}
\hline
$\bolds{\nu}$ & $\bolds{E\llvert \hat{\nu}_{b,2} - \nu\rrvert}$ & $\bolds{E\llvert \hat{\nu}_b - \nu\rrvert}$ \\
\hline
$0.1$ & 0.102 (0.007) & 0.059 (0.004) \\
$0.3$ & 0.133 (0.010) & 0.058 (0.006) \\
$0.5$ & 0.122 (0.009) & 0.058 (0.007) \\
$ 0.7$ & 0.111 (0.008) & 0.082 (0.007) \\
$0.9$ & 0.078 (0.007) & 0.068 (0.005) \\
$1.0$ & 0.076 (0.007) & 0.070 (0.005) \\
$1.3$ & 0.059 (0.004) & 0.059 (0.004) \\
$1.5$ & 0.049 (0.004) & 0.049 (0.004) \\
$1.7$ & 0.041 (0.003) & 0.041 (0.003) \\
$1.9$ & 0.072 (0.004) & 0.072 (0.004) \\ \hline
\end{tabular*}
\end{table}

\item[(ii)] $\gamma$ is an arc of the unit circle given by $\gamma(t) = (\cos
(t), \sin(t) )'$ for all $0\leq t \leq\pi/2$.
\end{longlist}
\end{rem}

\begin{expp}\label{expp2}
Set $n = 200$, $\nu= 0.1$, 0.3, 0.5, 0.7, 0.9, 1,
1.3, 1.5, 1.7, 1.9 and let $\varphi(s) = s(s+1)/(L+1)$,
$0\leq s\leq L=\pi/2$.
The data $X_1,\ldots, X_n$ are given by $X_i = X( \gamma( t_i ) )$
where $t_i = \varphi( L (i-1)/(n-1) )$, $i=1,\ldots, n$.
Table~\ref{tab2} reports the estimated mean absolute errors of $\hat{\nu}_{b,2}$
and $\hat{\nu}_b$ over 100 replications.
The mean absolute errors of $\hat{\nu}_{b,2}$ and
$\hat{\nu}_b$ indicate that $\hat{\nu}_b$ is indeed more accurate than
$\hat{\nu}_{b,2}$.
\end{expp}

\section{Second-order quadratic variations using $4$ or more deformed lattice points in ${\mathbb R}^2$}\label{sec4}

Let $[0,1]^2 \subset\Omega$ where $\Omega$ is an open set in ${\mathbb R}^2$.
A function $\tilde{\varphi}: \Omega\rightarrow {\mathbb R}^2$ is said
to be $C^1 (\Omega)$ if all first-order partial derivatives $\partial
\varphi_i (x_1, x_2)/\partial x_j$ exist and are continuous.
Here, we write $\tilde{\varphi} = (\varphi_1, \varphi_2)$ where $\varphi
_1, \varphi_2$ are real-valued functions.
The class of $C^1 (\Omega)$ diffeomorphisms is the set of all
continuous invertible maps
$\tilde{\varphi}: \Omega\rightarrow {\mathbb R}^2$, such that
$\tilde{\varphi}$ is $C^1(\Omega)$ and $\tilde{\varphi}^{-1}$ is $C^1
(\Omega^*)$ where $\Omega^*$ is the range of $\tilde{\varphi}$.

\begin{co}\label{co5}
$\tilde{\varphi}$ is a $C^2 (\Omega)$
diffeomorphism, which is a $C^1 (\Omega)$ diffeomorphism with
continuous second-order mixed partial derivatives.

Define ${\mathbf x}^{i_1, i_2} = (\varphi_1 (i_1/n, i_2/n), \varphi_2
(i_1/n, i_2/n))'$, $1 \leq i_1, i_2 \leq n$,
and denote $X_{i_1, i_2} = X({\mathbf x}^{i_1, i_2})$ where $X$ is a
stationary, isotropic Gaussian random field on ${\mathbb R}^2$ having
covariance function $K$ as in (\ref{eq:1.01}) with $d=2$.
We write $\varphi_j^{(1,0)} (u, v) = \partial\varphi_j (u, v)/\partial
u$ and
$\varphi_j^{(0,1)} (u, v) = \partial\varphi_j (u, v)/\partial v$ for $j=1,2$.

Since $\tilde{\varphi}$ is a $C^2(\Omega)$ diffeomorphism, we observe
from \cite{ac1}, page 2331, that
there exist constants $0 < C_{5, 0} \leq1 \leq C_{5, 1}$ such that
that for all $1\leq i_1, i_2,j_1, j_2 \leq n$,
\[
C_{5, 0} \biggl\llVert \biggl( \frac{i_1}{n}, \frac{ i_2}{n}
\biggr)' - \biggl(\frac{j_1}{n}, \frac{j_2}{n}
\biggr)' \biggr\rrVert\leq \bigl\llVert{\mathbf x}^{i_1, i_2} -
{ \mathbf x}^{ j_1, j_2} \bigr\rrVert\leq C_{5, 1} \biggl\llVert
\biggl( \frac{ i_1}{n}, \frac{i_2}{n} \biggr)' - \biggl(
\frac{j_1}{n}, \frac{j_2}{n} \biggr)' \biggr\rrVert.
\]

For $\theta\in\{1,2\}$ and $1\leq i_1, i_2 \leq n- \theta$, let
\begin{eqnarray*}
{\mathbf x}^{i_1+ \theta k_1, i_2 + \theta k_2} &=& \bigl( x_1^{i_1+
\theta k_1, i_2 + \theta k_2},
x_2^{i_1+ \theta k_1, i_2 + \theta k_2} \bigr)' \qquad\forall0 \leq
k_1, k_2 \leq1,
\nonumber
\\
A_{\theta; i_1, i_2} &=& \pmatrix{ x_1^{i_1+\theta, i_2}-
x_1^{i_1, i_2} & x_2^{i_1+\theta, i_2} -
x_2^{i_1, i_2} \vspace*{3pt}
\cr
x_1^{i_1, i_2+ \theta}-
x_1^{i_1, i_2} & x_2^{i_1, i_2+ \theta} -
x_2^{i_1, i_2} },
\nonumber
\\
B_{\theta; i_1, i_2} &=& \pmatrix{ x_1^{i_1+ \theta, i_2}-
x_1^{i_1+ \theta, i_2+\theta} & x_2^{i_1+ \theta, i_2} -
x_2^{i_1+ \theta, i_2+ \theta} \vspace*{3pt}
\cr
x_1^{i_1, i_2+ \theta}-
x_1^{i_1+ \theta, i_2+ \theta} & x_2^{i_1, i_2+ \theta} -
x_2^{i_1+ \theta, i_2+\theta} }.
\end{eqnarray*}
We write
%
%
\begin{equation}
\pmatrix{ X_{i_1+\theta, i_2} - X_{i_1, i_2} \vspace*{3pt}
\cr
X_{i_1, i_2+\theta} - X_{i_1, i_2} } = A_{\theta; i_1, i_2} \pmatrix{
g_{\theta, 1} (i_1, i_2) \vspace*{3pt}
\cr
g_{\theta, 2} (i_1, i_2) }. \label{eq:8.1}
\end{equation}
The motivation for (\ref{eq:8.1}) is that $g_{\theta, j} (i_1, i_2)$
approximates $\partial X_{i_1, i_2}/\partial x_j^{i_1, i_2}$.
Since $\tilde{\varphi}$ is a diffeomorphism, it follows from the
inverse function theorem (cf. \cite{l2}) that
\begin{eqnarray*}
(n/\theta)^2 \llvert A_{\theta; i_1, i_2} \rrvert&\rightarrow&
\varphi_1^{(1,0)} \biggl(\frac{i_1}{n}, \frac{i_2}{n}
\biggr) \varphi_2^{(0,1)} \biggl(\frac{i_1}{n},
\frac{i_2}{n} \biggr)
- \varphi_1^{(0,1)} \biggl( \frac{i_1}{n},
\frac{i_2}{n} \biggr) \varphi_2^{(1, 0)} \biggl(
\frac{i_1}{n}, \frac{i_2}{n} \biggr)
\\
&\neq& 0,
\end{eqnarray*}
as $n\rightarrow\infty$. Thus, for sufficiently large $n$, $A_{\theta;
i_1, i_2}^{-1}$ exists and
\begin{eqnarray*}
A_{\theta; i_1, i_2}^{-1} &=& \llvert A_{\theta; i_1,i_2}\rrvert
^{-1} \pmatrix{ x_2^{i_1, i_2+ \theta}- x_2^{i_1, i_2}
& -x_2^{i_1+ \theta, i_2} + x_2^{i_1, i_2} \vspace*{3pt}
\cr
-x_1^{i_1, i_2+ \theta} + x_1^{i_1, i_2} &
x_1^{i_1+ \theta, i_2} - x_1^{i_1, i_2} }
\nonumber
\\
&=& \bigl( \alpha_{\theta; i_1, i_2}^{j,l} \bigr)_{j, l =1, 2}, \qquad
\mbox{say}.
\end{eqnarray*}
Hence,
\[
\pmatrix{ g_{\theta, 1} (i_1, i_2) \vspace*{3pt}
\cr
g_{\theta, 2} (i_1, i_2) } = A_{\theta; i_1, i_2}^{-1}
\pmatrix{ X_{i_1+ \theta, i_2} - X_{i_1, i_2} \vspace*{3pt}
\cr
X_{i_1, i_2+ \theta} - X_{i_1, i_2} }.
\]
In particular choosing the $\ell$th coordinate, we obtain
%
\begin{eqnarray}
g_{\theta, \ell} (i_1, i_2) = \alpha_{\theta; i_1, i_2}^{\ell, 1}
( X_{i_1 + \theta, i_2} - X_{i_1, i_2} ) + \alpha_{\theta; i_1,
i_2}^{\ell, 2}
( X_{i_1, i_2+ \theta} - X_{i_1, i_2} ),
\nonumber
\\
\eqntext{\ell\in\{1,2\}.}
\end{eqnarray}
In a similar manner, writing
%
%
\begin{equation}
\pmatrix{ X_{i_1+ \theta, i_2} - X_{i_1+ \theta, i_2+ \theta} \vspace*{3pt}
\cr
X_{i_1, i_2+ \theta} - X_{i_1+ \theta, i_2+ \theta} } = B_{\theta; i_1, i_2} \pmatrix{
h_{\theta, 1} (i_1+ \theta, i_2+ \theta)
\vspace*{3pt}
\cr
h_{\theta, 2} (i_1+ \theta, i_2+
\theta) }, \label{eq:8.2}
\end{equation}
we have
\begin{eqnarray*}
B_{\theta; i_1, i_2}^{-1} &=& \llvert B_{\theta; i_1,i_2}\rrvert
^{-1} \pmatrix{ x_2^{i_1, i_2+ \theta}- x_2^{i_1+ \theta, i_2+ \theta}
& -x_2^{i_1+ \theta, i_2} + x_2^{i_1+ \theta, i_2+ \theta} \vspace*{3pt}
\cr
-x_1^{i_1, i_2+ \theta} + x_1^{i_1+ \theta, i_2+ \theta} &
x_1^{i_1+ \theta, i_2} - x_1^{i_1+ \theta, i_2+ \theta} }
\nonumber
\\
&=& \bigl( \beta_{\theta; i_1, i_2}^{j,l} \bigr)_{j, l =1, 2}, \qquad
\mbox{say},
\end{eqnarray*}
and hence
\begin{eqnarray*}
\pmatrix{ h_{\theta, 1} (i_1+ \theta, i_2+ \theta)
\vspace*{3pt}
\cr
h_{\theta, 2} (i_1+ \theta, i_2+
\theta) } & = & B_{\theta; i_1, i_2}^{-1} \pmatrix{ X_{i_1+ \theta, i_2} -
X_{i_1+
\theta, i_2+ \theta} \vspace*{3pt}
\cr
X_{i_1, i_2+ \theta} - X_{i_1+ \theta, i_2+ \theta} },
\nonumber
\\
h_{\theta, \ell} (i_1+\theta, i_2+\theta) &=&
\beta_{\theta; i_1, i_2}^{\ell, 1} ( X_{i_1+ \theta, i_2 } - X_{i_1+
\theta, i_2+ \theta} )
\nonumber
\\
&&{} + \beta_{\theta; i_1, i_2}^{\ell, 2} ( X_{i_1, i_2+ \theta} -
X_{i_1+ \theta, i_2+ \theta} )\qquad\forall\ell\in\{1,2\}.
\end{eqnarray*}
For $1\leq i_1, i_2 \leq n-\theta$ and $\theta, \ell\in\{1,2\}$, we write
%
%
\begin{eqnarray}
\label{eq:8.3} \bar{\nabla}_{\theta, \ell} X_{i_1, i_2} &=&
\beta_{\theta; i_1,
i_2}^{\ell, 1} ( X_{i_1+ \theta, i_2} - X_{i_1+ \theta, i_2+ \theta} ) +
\beta_{\theta
; i_1, i_2}^{\ell, 2} ( X_{i_1, i_2+ \theta} - X_{i_1+ \theta, i_2+ \theta} )
\nonumber
\\
&&{} - \alpha_{\theta; i_1, i_2}^{\ell, 1} ( X_{i_1+ \theta, i_2} -
X_{ i_1, i_2} ) - \alpha_{\theta; i_1, i_2}^{\ell, 2} ( X_{i_1, i_2+
\theta} -
X_{ i_1, i_2} )
\\
&=& \sum_{0\leq k_1, k_2 \leq1} c_{\theta, \ell; i_1, i_2}^{k_1, k_2}
X_{i_1 + \theta k_1, i_2+ \theta k_2}, \qquad\mbox{say}.
\nonumber
\end{eqnarray}
In this section, second-order quadratic variations based on $\{ X_{i_1,
i_2}: 1\leq i_1, i_2 \leq n\}$ are defined to be
%
%
\begin{equation}
\bar{V}_{\theta, \ell} = \sum_{1\leq i_1, i_2 \leq n-\theta} ( \bar{
\nabla}_{\theta, \ell} X_{i_1, i_2} )^2, \qquad\theta, \ell\in
\{1,2 \}. \label{eq:7.2}
\end{equation}

Lemma~\ref{la:7.1} below provides the rationale for calling $\bar
{V}_{\theta, \ell}$ ``second-order''.
\end{co}

%
\begin{la} \label{la:7.1}
Let $c_{\theta, \ell; i_1, i_2}^{k_1, k_2}$, $\theta, \ell\in\{1,2 \}
$, be as in (\ref{eq:8.3}).
Then for all $1\leq i_1, i_2 \leq n-\theta$, we have
\begin{eqnarray*}
\sum_{0\leq k_1, k_2 \leq1} c_{\theta, \ell; i_1, i_2}^{k_1, k_2} &=& 0,
\nonumber
\\
\sum_{0\leq k_1, k_2 \leq1} c_{\theta, \ell; i_1, i_2}^{k_1, k_2}
x_j^{i_1+ \theta k_1, i_2+ \theta k_2} &=& 0 \qquad\forall j = 1, 2.
\end{eqnarray*}
\end{la}

Let $G_\nu(\cdot)$ as in (\ref{eq:1.19}). For $\theta, \ell\in\{1,2\}$, define
%
%
\begin{eqnarray}
\label{eq:4.00} \bar{f}_{\theta, \ell} (\nu) &=& \beta_\nu^* \sum
_{1\leq i_1, i_2\leq n-\theta} \sum_{0\leq k_1, l_1, k_2, l_2 \leq1:
(k_1,k_2)\neq(l_1, l_2)}
c_{\theta, \ell; i_1, i_2}^{k_1, k_2} c_{\theta, \ell; i_1, i_2}^{l_1, l_2}
\nonumber
\\[-8pt]
\\[-8pt]
\nonumber
&&{}\times G_\nu \bigl( \bigl\llVert{\mathbf
x}^{i_1 + \theta k_1, i_2 +
\theta k_2} - {\mathbf x}^{i_1 + \theta l_1, i_2 + \theta l_2} \bigr\rrVert \bigr)\qquad \forall
\nu\in(0,2),
\end{eqnarray}
and $\tilde{K} ({\mathbf x}- {\mathbf y}) = K( {\mathbf x}, {\mathbf y})$ for all ${\mathbf
x}, {\mathbf y} \in{\mathbb R}^2$. Letting $a_1, a_2$ to be nonnegative integers,
we write
\[
\tilde{K}^{(a_1, a_2)} ({\mathbf x}) = \frac{\partial^{a_1 + a_2} }{
\partial x_1^{a_1}\,\partial x_2^{a_2} } \tilde{K} ({\mathbf
x}) \qquad\forall{\mathbf x} = (x_1, x_2)',
\]
and $\Delta= \{ (\varphi_1 (u,v), \varphi_2 (u,v) )': 0\leq u,v\leq1\}$.

%
\begin{tm} \label{tm:7.1}
Let $\theta, \ell\in\{1, 2\}$, $\nu\in(0, 2)$ and $\bar{V}_{\theta,
\ell}$ be as in (\ref{eq:7.2}).
Suppose that
Condition~\ref{co5} holds and that ${\mathcal J}_{\nu, \ell} (\cdot,\cdot)$, as in
Lemma~\ref{la:7.11}, is not identically $0$ on $[0,1]^2$.
Then
$E( \bar{V}_{\theta, \ell} ) \sim\bar{f}_{\theta, \ell} (\nu) \asymp
n^{4 - 2\nu}$.
In addition, suppose \textup{(i)}~$\tilde{K}^{(a_1, a_2)} (\cdot)$ is a continuous
function on an open set containing $\{ {\mathbf x}-{\mathbf y}: {\mathbf x}, {\mathbf
y}\in\Delta, {\mathbf x}\neq{\mathbf y}\}$
whenever $a_1, a_2$ are nonnegative integers such that $a_1 + a_2= 4$ and
\textup{(ii)}~there exists a constant $\tilde{C}$ such that
\[
\bigl\llvert\tilde{K}^{(a_1, a_2)} ( {\mathbf x} - {\mathbf y}) \bigr\rrvert
\leq \tilde{C} \llVert{\mathbf x} - {\mathbf y} \rrVert^{2 \nu- 4},
\]
for all $a_1 + a_2 = 4$, ${\mathbf x}, {\mathbf y}\in \Delta$ and ${\mathbf x}\neq
{\mathbf y}$.
Then
\[
\operatorname{Var} \bigl\{ \bar{V}_{\theta, \ell}/ E (\bar{V}_{\theta, \ell} )
\bigr\} = \cases{ O \bigl( n^{ -2} \bigr), &\quad if $\nu\in(0,3/2)$,
\vspace*{3pt}
\cr
O \bigl\{ n^{-2} \log(n) \bigr\}, &\quad if $\nu= 3/2$,
\vspace*{3pt}
\cr
O \bigl( n^{- 8 +4 \nu} \bigr), &\quad if $\nu\in(3/2, 2)$,}
\]
and $\bar{V}_{\theta, \ell}/E( \bar{V}_{\theta, \ell} ) \rightarrow1$
almost surely
as $n\rightarrow\infty$.
\end{tm}

The proof of Theorem~\ref{tm:7.1}, which uses Lemmas~\ref{la:7.11} and
\ref{la:7.13}, is similar to that of Theorem~\ref{tm:4.1} and can be
found in \cite{l2}.
We end this section with the following immediate corollary of Theorem
\ref{tm:7.1}.

\begin{cy} \label{cy:7.1}
Suppose $\nu\in(0,2)$ and $\theta, \ell\in\{1, 2\}$. Then under the
conditions of Theorem~\ref{tm:7.1}, $\hat{\nu}_{\theta, \ell} = \{ 4
-\log ( \bar{V}_{\theta, \ell})/\log(n) \}/2$ is a strongly
consistent estimator for $\nu$.
\end{cy}

We observe that, unfortunately, the bias of $\hat{\nu}_{\theta, \ell}$
is of order $1/\log(n)$ and this makes $\hat{\nu}_{\theta, \ell}$
unsuitable for use in practice.

\subsection{Estimating \texorpdfstring{$\nu$}{$nu$} using deformed lattice data in ${\mathbb R}^2$}\label{sec4.1}

Writing $\bar{f}_{\theta, \ell} (\cdot)$, $\theta, \ell\in\{1,2\}$, as in
(\ref{eq:4.00}),
define $\bar{F}_{\ell, n}: [0, 2] \rightarrow[0, \infty)$ by
\[
\bar{F}_{\ell, n} \bigl(\nu^*\bigr) = \bar{f}_{2, \ell} \bigl(\nu^*
\bigr)/ \bar{f}_{1, \ell} \bigl(\nu^*\bigr)\qquad\mbox{if $\nu ^*
\in(0,2)$},
\]
$\bar{F}_{\ell, n} (0) = \lim_{\delta\rightarrow0+} \bar{F}_{\ell, n}
(\delta)$ and $\bar{F}_{\ell, n} (2) = \lim_{\delta\rightarrow0+} \bar
{F}_{\ell, n} (2 -\delta)$.
We observe from Lemma~\ref{la:a.4} that $\bar{F}_{\ell, n} (\cdot)$ is a
continuous function on $[0,2]$.
Let $\hat{\nu}_{c, \ell} \in[0,2]$ satisfy
%
%
\begin{equation}
\biggl\{ \frac{ \bar{V}_{1, \ell} \bar{F}_{\ell, n} (\hat{\nu}_{c, \ell
}) }{ \bar{V}_{2, \ell} } - 1 \biggr\}^2 = \min
_{0 \leq\nu^* \leq2} \biggl\{ \frac{ \bar{V}_{1, \ell} \bar{F}_{\ell,
n} (\nu^*) }{ \bar{V}_{2, \ell}} - 1 \biggr\}^2.
\label{eq:10.1}
\end{equation}

The following is the main result of this section.

\begin{tm} \label{tm:10.1}
Let $\nu\in(0, 2)$, $\ell\in\{1,2\}$ and $\hat{\nu}_{c, \ell}$ be
as in (\ref{eq:10.1}). Suppose
the conditions of Theorem~\ref{tm:7.1} are satisfied.
Then $\hat{\nu}_{c, \ell} \rightarrow\nu$ as $n\rightarrow\infty$
almost surely.
\end{tm}

Using Lemma~\ref{la:7.13}, the proof of Theorem~\ref{tm:10.1} is
similar to that of Theorem~\ref{tm:2.2} and can be found in \cite{l2}.

\begin{rem}\label{rem5}
We observe from the definition of $\hat{\nu}_{c, \ell}$
that $\tilde{\varphi}$ need not be explicitly known; only the bijection
$\tilde{\Phi}: \{(i_1, i_2): 1 \leq i_1, i_2\leq n\} \rightarrow\{\phi
_{i_1, i_2} = \tilde{\varphi} (i_1/n, i_2/n): 1\leq i_1, i_2 \leq n\}$
is required.
In many designs of experiments scenarios, $\tilde{\varphi}$ is known.
On the other hand, if $\tilde{\varphi}$ is unknown, we propose the
algorithm below to determine the bijection $\tilde{\Phi}$.

\begin{longlist}[\textit{Step} 4.1.]
\item[\textit{Step} 4.1.] First, divide the set $\{ \phi_{i_1,i_2}, 1 \leq i_1,
i_2 \leq n\}$ into $n$ disjoint subsets $\gamma_1, \ldots, \gamma_n$
where each subset contains $n$ elements.
The motivation here is that $\gamma_i$ corresponds to $\{ \phi_{i,
j}:1\leq j\leq n\}$.
As in Remark~\ref{rem4}, the elements of $\gamma_i$ are ordered and adjacent
elements are connected by straight lines. Thus,
$\gamma_i$ becomes a piecewise linear continuous curve in ${\mathbb
R}^2$. As curves, $\gamma_1, \ldots, \gamma_n$ are chosen so that
$\gamma_i$ and $\gamma_j$ do not intersect if $i\neq j$ and for each
$1\leq i\leq n-2$,
the only curve that lies between $\gamma_i$ and $\gamma_{i+2}$ is
$\gamma_{i+1}$.

\item[\textit{Step} 4.2.] For $\gamma_1$, we label its elements by ${\mathbf
x}^{1,1}, \ldots, {\mathbf x}^{1,n}$ in that the ordering of the second
superscript follows the ordering of the elements along the curve $\gamma_1$.

\item[\textit{Step} 4.3.] Let $1\leq k\leq n-1$. Given $\gamma_l = \{ {\mathbf x}^{l,
j}: 1\leq j\leq n\}$, $1\leq l\leq k$, the elements of $\gamma_{k+1}$
are ordered as in Remark~\ref{rem4} such that ${\mathbf x}^{k+1,1}$ is closer to
${\mathbf x}^{k, 1}$ than to~${\mathbf x}^{k, n}$.

\item[\textit{Step} 4.4.] Increase $k$ by 1 and repeat step 4.3 until the
elements of $\gamma_n$ have been ordered. The required bijection $\tilde
{\Phi}$ is given as ${\mathbf x}^{i,j} = \tilde{\varphi} (i/n, j/n)$,
$1\leq i, j \leq n$.
\end{longlist}

We note that there is an extensive literature on landmark matching via
large deformation diffeomorphisms, for example, \cite{gm1,jm1}.
These techniques can be used to compute the diffeomorphism $\tilde
{\varphi}$ in step 4.4 with bijection $\tilde{\Phi}$ as landmarks.
\end{rem}

\begin{rem}\label{rem6}
The algorithm in Remark~\ref{rem5} assumes that the curves
$\gamma_1, \ldots,\break \gamma_n$ can be found. A simpler alternative is to
select only 1 subset $\Gamma\subseteq\{\tilde{\varphi} (i_1/n,\break
i_2/n): 1\leq i_1, i_2 \leq n\}$
of cardinality $n$, say. Treat $\Gamma$ as points lying on a $C^2$
curve $\gamma$ and then apply the methodology of Section~\ref{sec3} to estimate
$\nu$.
For the asymptotics of Theorem~\ref{tm:9.10} to be more effective,
$\Gamma$ should be chosen so that the adjacent points on $\gamma$ are
close to each other and that the curvature of the curve $\gamma$
is reasonably small.

A simulation experiment is conducted to gauge the finite sample
accuracy of $\hat{\nu}_{c, 1}$ and $\hat{\nu}_{c, 2}$.
In Experiment~\ref{expp3} below,
$X$ is a stationary Gaussian random field having mean 0 and Mat\'{e}rn
covariance function $K_{\mathrm{Mat}}$ as in (\ref{eq:2.1.1}) with $\sigma=
\alpha= 1$ and $d = 2$.
\end{rem}

\begin{expp}\label{expp3}
Set $n = 40$ and $\nu= 0.1$, 0.3, 0.5, 0.7, 0.9,
1, 1.3, 1.5, 1.7, 1.9.
Let $\tilde{\varphi}: {\mathbb C}\rightarrow{\mathbb C}$ be given by
$\tilde{\varphi} (z) = z (z+1)/3 = \varphi_1 (z) + {\mathbf i} \varphi
_2(z)$ $\forall z \in{\mathbb C}$
where $\varphi_i$'s are real-valued functions and ${\mathbf i}=\sqrt{-1}$.
We take ${\mathbf x}^{i_1,i_2} = (\varphi_1 (i_1 n^{-1}+ {\mathbf i} i_2
n^{-1}), \varphi_2 (i_1 n^{-1} + {\mathbf i} i_2 n^{-1}))'$, $1\leq i_1,
i_2\leq n$.
Table~\ref{tab3} reports the estimated mean absolute errors of $\hat{\nu}_{c,
1}$ and $\hat{\nu}_{c, 2}$
over 100 replications.
\end{expp}

%
\begin{table}
\tabcolsep=0pt
\tablewidth=233pt
\caption{Experiment \protect\ref{expp3}. Simulation results in estimating $\nu$
using deformed lattice data over 100 replications (standard errors
within parentheses)}\label{tab3}
\begin{tabular*}{\tablewidth}{@{\extracolsep{\fill}}@{}lcc@{}}
\hline
$\bolds{\nu}$ & $\bolds{E\llvert \hat{\nu}_{c, 1} - \nu\rrvert}$ & $\bolds{E \llvert \hat {\nu}_{c, 2} - \nu\rrvert}$ \\
\hline
$0.1$ & 0.054 (0.003) & 0.036 (0.003) \\
$0.3$ & 0.042 (0.004) & 0.034 (0.002) \\
$0.5$ & 0.051 (0.003) & 0.034 (0.003) \\
$ 0.7$ & 0.041 (0.003) & 0.033 (0.002) \\
$0.9$ & 0.040 (0.003) & 0.029 (0.002) \\
$1.0$ & 0.044 (0.003) & 0.033 (0.002) \\
$1.3$ & 0.040 (0.003) & 0.030 (0.002) \\
$1.5$ & 0.036 (0.003) & 0.034 (0.003) \\
$1.7$ & 0.064 (0.005) & 0.056 (0.004) \\
$1.9$ & 0.118 (0.005) & 0.102 (0.004) \\
\hline
\end{tabular*}
\end{table}

\begin{appendix}
\section*{Appendix}\label{appen}

%
\begin{la} \label{la:a.1}
Let $\theta\in\{1,2\}$ and $f_{\theta, \ell}(\cdot)$, $F_{\ell, n}(\cdot)$ be
as in (\ref{eq:2.101}), (\ref{eq:2.1.4}), respectively.
Suppose $0 < \nu\leq M < \ell\leq10$ and Condition~\ref{co1} holds. Then
\begin{eqnarray*}
f_{\theta, \ell} ( \nu) &= & 2 \beta_\nu^* \theta^{2 \nu-2 \ell}
n^{2 \ell+1 - 2\nu} H_\ell(\nu) \int_0^1
\bigl\{ \varphi^{(1)} ( s) \bigr\}^{2 \nu-2 \ell} \,ds + O
\bigl(n^{2 \ell- 2\nu} \bigr)
\nonumber
\\
& \asymp& n^{2\ell+ 1 - 2 \nu},
\end{eqnarray*}
and hence $F_{\ell, n} (\nu^*) = 2^{2 \nu^* - 2 \ell} + O(n^{-1})$
as $n\rightarrow\infty$ uniformly over $0 \leq\nu^* \leq M$ where
\[
H_\ell(\nu) = \sum_{0\leq k_1 <k_2\leq\ell}
(-1)^{k_1 + k_2} \pmatrix{\ell
\cr
k_1} \pmatrix{\ell
\cr
k_2} G_\nu( k_2 - k_1 ) \qquad
\forall\nu\in (0, M].
\]
\end{la}

\begin{pf}
Since $\varphi$ is twice continuously differentiable, we
have $t_{i+ \theta k_2} - t_{i+ \theta k_1} = \theta(k_2 - k_1)
(n-1)^{-1} \varphi^{(1)} ((i-1)/(n-1)) + O(n^{-2})$ as $n\rightarrow
\infty$ uniformly over $1\leq i \leq n- \theta\ell$.
It follows from (\ref{eq:2.101}) that
\begin{eqnarray*}
f_{\theta, \ell} (\nu) &= & 2 \beta_\nu^* \theta^{2 \nu- 2 \ell}
n^{2 \ell - 2\nu} H_\ell(\nu) \sum_{i=1}^{n- \theta\ell}
\biggl\{ \varphi^{(1)} \biggl( \frac{i-1 }{n-1} \biggr) \biggr
\}^{2 \nu-2 \ell} + O \bigl(n^{2 \ell- 2 \nu} \bigr)
\nonumber
\\
&= & 2 \beta_\nu^* \theta^{2 \nu-2 \ell} n^{2 \ell+1 - 2\nu}
H_\ell(\nu) \int_0^1 \bigl\{
\varphi^{(1)} ( s) \bigr\}^{2 \nu-2 \ell} \,ds + O \bigl(n^{2 \ell- 2 \nu}
\bigr),
\end{eqnarray*}
as $n\rightarrow\infty$ uniformly over $0 < \nu\leq M$.
For each $\ell= 1,\ldots, 10$, we use {\em Mathematica} to plot
$H_\ell(\nu)$
to verify that $H_\ell(\nu) \neq0$ for all $0\leq\nu\leq M$.
\end{pf}

%
\begin{la} \label{la:a.2}
Let $f:{\mathbb R} \rightarrow{\mathbb R}$ and that for some
nonnegative integer $m$, $f^{(m+1)}$ is continuous on an open interval
containing $a$ and $x$. Then
\[
f(x) = \sum_{j=0}^m \frac{ f^{(j)} (a)}{ j!}
(x-a)^j + \frac{1}{m!} \int_a^x
(x-t)^m f^{(m+1)} (t) \,dt.
\]
\end{la}

\begin{pf*}{Proof of Theorem~\ref{tm:2.1}}
Without loss of generality, we shall assume that $E(X_i) = 0$.

\begin{longlist}[(a)]
\item[(a)] We observe from (\ref{eq:1.01}) that
\[
K( x, y) = \sum_{j = 0}^{ \lfloor\nu\rfloor}
\beta_j (x-y)^{2 j} + \beta_\nu^*
G_\nu \bigl(\llvert x-y \rrvert \bigr) + r(x, y)\qquad\forall x,y \in
\mathbb R.
\]
Since $\nu\in(0, \ell)$, $f_{\theta, \ell} (\nu) \asymp n^{2 \ell+1
- 2 \nu}$ (cf. Lemma~\ref{la:a.1}). As
$r( x, y ) = O( \llvert x -y \rrvert ^{2 \nu+ \tau} )$ for some
constant $\tau>0$ as $\llvert x-y\rrvert \rightarrow0$, it follows
from Lemma~\ref{la:a.01} that
\[
\sum_{i=1}^{n- \theta\ell} E \bigl\{ (
\nabla_{\theta, \ell} X_i)^2 \bigr\} = \sum
_{i=1}^{n- \theta\ell} \sum_{k_1 =0}^\ell
\sum_{k_2 =0}^\ell a_{\theta, \ell; i, k_1}
a_{\theta, \ell; i, k_2 } K( t_{i+ \theta k_1}, t_{i+ \theta k_2} ) \sim
f_{\theta, \ell} ( \nu).
\]
Hence, $E(V_{\theta,\ell}) \sim f_{\theta, \ell} (\nu)$.
Similarly, we observe from (\ref{eq:2.4}) and Lemma~\ref{la:a.01} that
\begin{eqnarray*}
&& \sum_{j= i - \theta\ell-1}^{i + \theta\ell+ 1} \llvert
\Sigma_{ i, j} \rrvert
\nonumber
\\
&&\qquad = \frac{1}{ E (V_{\theta, \ell} )} \sum_{j= i - \theta\ell-1}^{i +
\theta\ell+ 1}
\Biggl\llvert\sum_{k_1 =0}^\ell\sum
_{k_2 =0}^\ell a_{\theta, \ell; i, k_1} a_{\theta, \ell; j, k_2 }
\tilde{K} ( t_{i+ \theta k_1} - t_{ j + \theta k_2 } ) \Biggr\rrvert
\nonumber
\\
&&\qquad = \frac{O(1) }{ E (V_{\theta, \ell} )} \sum_{j= i - \theta\ell
-1}^{i + \theta\ell+ 1}
\Biggl\llvert\sum_{k_1 =0}^\ell\sum
_{k_2 =0}^\ell a_{\theta, \ell; i, k_1} a_{\theta, \ell; j, k_2 }
G_\nu \bigl(\llvert t_{i + \theta k_1} - t_{ j + \theta k_2 }\rrvert
\bigr) \Biggr\rrvert
\\
&&\qquad = O \bigl( n^{-1} \bigr),
\end{eqnarray*}
as $n\rightarrow\infty$ uniformly over $1\leq i \leq n- \theta\ell$.
Here, for brevity, we write\break $\sum_{j= i - \theta\ell-1}^{i + \theta
\ell+ 1} = \sum_{j= (i - \theta\ell-1)\vee1}^{(i + \theta\ell+
1)\wedge(n- \theta\ell)}$
since $j\in\{1,\ldots, n- \theta\ell\}$.
For $1\leq i, j\leq n - \theta \ell$, we observe from Lemmas~\ref
{la:a.01} and~\ref{la:a.2} that
\begin{eqnarray*}
\Sigma_{i, j} &=& \frac{1}{ E (V_{\theta, \ell} )} E \Biggl\{ \Biggl( \sum
_{k_1 =0}^\ell a_{\theta, \ell; i, k_1} X_{i+ \theta k_1}
\Biggr) \Biggl( \sum_{k_2 =0}^\ell
a_{\theta, \ell; j, k_2} X_{j + \theta k_2} \Biggr) \Biggr\}
\nonumber
\\
&=& \frac{1}{ E (V_{\theta, \ell} )} \sum_{k_1 =0}^\ell
\sum_{k_2 =0}^\ell a_{\theta, \ell; i, k_1}
a_{\theta, \ell; j, k_2} \tilde{K} ( t_{j + \theta k_2} - t_{i+ \theta
k_1} )
\nonumber
\\
&=& \frac{1}{ (2 \ell-1)! E (V_{\theta, \ell} )} \sum_{k_1 =0}^\ell
\sum_{k_2 =0}^\ell a_{\theta, \ell; i, k_1}
a_{\theta, \ell; j, k_2}
\nonumber
\\
&&{}\times\int_{t_j - t_i}^{t_{j+ \theta k_2} - t_{i+ \theta k_1}} (t_{j+\theta k_2} -
t_{i+ \theta k_1} -t)^{2 \ell-1} \tilde{K}^{(2 \ell)} (t) \,dt,
\end{eqnarray*}
if the right-hand side exists. Since $\llvert K^{( 2 \ell)}
(t)\rrvert \leq C_\ell t^{2\nu- 2 \ell}$, $t\in(0,1]$, it follows
from Condition~\ref{co1} that
\begin{eqnarray*}
\sum_{j= 1}^{i- \theta\ell-2} \llvert
\Sigma_{i, j} \rrvert&= & \frac{1}{ (2 \ell-1)! E (V_{\theta, \ell}
)} \sum
_{j = 1}^{i- \theta\ell-2} \Biggl\llvert\sum
_{k_1 =0}^\ell\sum_{k_2 =0}^\ell
a_{\theta, \ell; i, k_1} a_{\theta, \ell; j, k_2 }
\\
&&{}\times \int_{t_i - t_j}^{
t_{i + \theta k_1} - t_{j + \theta k_2} } ( t_{i+ \theta k_1} -
t_{j + \theta k_2} - t)^{2 \ell-1} \tilde{K}^{( 2 \ell)} ( t) \,d t \Biggr
\rrvert
\nonumber
\\
&\leq& \frac{O(1)}{ E (V_{\theta, \ell} ) } \sum_{j=1}^{i- \theta\ell-2}
\sum_{k_1 =0}^\ell\sum
_{k_2 =0}^\ell \bigl\{ (t_{i+ \theta k_1} -
t_{j+ \theta k_2}) \wedge(t_i - t_j) \bigr
\}^{2 \nu- 2 \ell}
\nonumber
\\
&\leq& \frac{O(1)}{ E (V_{\theta, \ell} ) } \sum_{j=1}^{i- \theta\ell
-2 }
\biggl(\frac{ i - j - \theta\ell}{n} \biggr)^{2 \nu- 2 \ell}
\nonumber
\\
&\leq& \frac{O(n )}{ E (V_{\theta, \ell} ) } \int_{1/n}^1
t^{2 \nu- 2 \ell} \,d t
\nonumber
\\
&=& \cases{ O \bigl( n^{ - 1} \bigr), &\quad if $\nu< (2 \ell-1)/2$,
\vspace*{3pt}
\cr
O \bigl\{ n^{ -1} \log(n) \bigr\}, &\quad if $\nu= (2
\ell-1)/2$, \vspace*{3pt}
\cr
O \bigl( n^{ - 2 \ell+ 2 \nu} \bigr), &\quad if $\nu> (2
\ell-1)/2$,} 
\end{eqnarray*}
as $n\rightarrow\infty$ uniformly over $\theta\ell+3 \leq i\leq n-
\theta\ell$. Similarly,
\[
\sum_{j= i+ \theta\ell+2 }^{n- \theta\ell} \llvert
\Sigma_{i, j} \rrvert= \cases{ O \bigl( n^{ - 1} \bigr), &\quad
if $\nu< (2 \ell-1)/2$, \vspace*{3pt}
\cr
O \bigl\{ n^{ -1} \log(n)
\bigr\}, &\quad if $\nu= (2 \ell-1)/2$, \vspace*{3pt}
\cr
O \bigl(
n^{ - 2 \ell+ 2 \nu} \bigr), &\quad if $\nu> (2 \ell-1)/2$,}
\]
as $n\rightarrow\infty$ uniformly over $1\leq i\leq n- 2 \theta\ell-2$.
Consequently, it follows from \cite{hj1} that
\begin{eqnarray}
\llVert\Sigma_{\mathrm{abs}} \rrVert_2 & \leq& \max
_{ 1 \leq i\leq n- \theta\ell} \Biggl\{\sum_{j=1}^{i- \theta\ell- 2}
\llvert\Sigma_{i, j}\rrvert+ \sum_{j= i- \theta\ell-1}^{i+ \theta
\ell+1}
\llvert\Sigma_{i,j}\rrvert+ \sum_{j=i+ \theta\ell+ 2}^{n- \theta
\ell}
\llvert\Sigma_{i, j}\rrvert \Biggr\}
\nonumber
\nonumber
\\[-8pt]
\\[-8pt]
\nonumber
& = & \cases{ O \bigl( n^{ - 1} \bigr), &\quad if $\nu< (2
\ell -1)/2$, \vspace*{3pt}
\cr
O \bigl\{ n^{ -1} \log(n) \bigr\}, &\quad
if $\nu= (2 \ell -1)/2$, \vspace*{3pt}
\cr
O \bigl( n^{ - 2 \ell+ 2 \nu} \bigr), &
\quad if $\nu> (2 \ell-1)/2$,} \label{eq:2.2}
\end{eqnarray}
as $n\rightarrow\infty$.
Next, we observe that
\[
\sum_{i=1}^{n- \theta\ell} \Sigma_{i,i}^2
+ \sum_{j=1}^{\theta\ell+1} \sum
_{i=1}^{n- \theta\ell-j} \Sigma_{i,i+j}^2 +
\sum_{j=1}^{\theta\ell+1} \sum
_{i= 1 +j}^{n- \theta\ell} \Sigma_{i,i-j}^2 = O
\bigl(n^{-1} \bigr),
\]
and
\begin{eqnarray*}
&& \sum_{i= \theta\ell+3}^{n- \theta\ell} \sum
_{j=1}^{i- \theta\ell-2} \Sigma_{i,j}^2
\nonumber
\\
&&\qquad \leq \frac{1}{ \{ (2\ell-1)! E (V_{\theta, \ell} ) \}^2} \sum_{i=
\theta\ell+3}^{n- \theta\ell}
\sum_{j=1}^{i- \theta\ell-2} \Biggl\{ \sum
_{k_1=0}^\ell\sum_{k_2 =0}^\ell
a_{\theta, \ell; i, k_1} a_{\theta, \ell; j, k_2}
\nonumber
\\
&&\quad\qquad{}\times\int_{t_i - t_j}^{ t_{i+ \theta k_1} - t_{j+ \theta k_2}} (
t_{i+ \theta k_1} - t_{j+ \theta k_2} -t)^{2 \ell-1} K^{(2 \ell)} ( t)
\,d t \Biggr\}^2
\nonumber
\\
&&\qquad \leq \frac{O(1)}{ \{ E (V_{\theta, \ell} ) \}^2 } \sum_{i= \theta
\ell+3}^{n- \theta\ell}
\sum_{j=1}^{i- \theta\ell-2} \Biggl[ \sum
_{k_1=0}^\ell\sum_{k_2=0}^\ell
\bigl\{ (t_{i+ \theta k_1} - t_{j+ \theta k_2}) \wedge( t_i -
t_j) \bigr\}^{2 \nu- 2 \ell} \Biggr]^2
\nonumber
\\
&&\qquad \leq \frac{O(1)}{ \{ E (V_{\theta, \ell}) \}^2 } \sum_{i= \theta
\ell+3}^{n- \theta\ell}
\sum_{j=1}^{i- \theta\ell-2} \biggl(\frac{ i-j- \theta\ell}{n}
\biggr)^{4 \nu-4 \ell}
\nonumber
\\
&&\qquad \leq \frac{O( n)}{ \{ E (V_{\theta, \ell} ) \}^2 } \sum_{i= \theta
\ell+3}^{n- \theta\ell}
\int_{1/n}^1 s^{4 \nu- 4 \ell} \,d s
\nonumber
\\
&&\qquad = \cases{ O \bigl( n^{-1} \bigr), &\quad if $\nu< (4
\ell-1)/4$, \vspace*{3pt}
\cr
O \bigl\{ n^{-1} \log(n) \bigr\}, &\quad if
$\nu= (4 \ell-1)/4$, \vspace*{3pt}
\cr
O \bigl( n^{-4 \ell+ 4 \nu} \bigr), &\quad if
$\nu> (4 \ell-1)/4$,}
\end{eqnarray*}
as $n\rightarrow\infty$. Hence, we conclude that
%
%
\begin{eqnarray}
\label{eq:2.8} \llVert\Sigma_{\mathrm{abs}}\rrVert^2_F
&=& \sum_{i=1}^{n- \theta\ell} \Sigma_{i,i}^2
+ \sum_{j=1}^{\theta\ell+1} \sum
_{i=1}^{n- \theta\ell-j} \Sigma_{i,i+j}^2 +
\sum_{j=1}^{\theta\ell+1} \sum
_{i= 1 +j}^{n- \theta\ell} \Sigma_{i,i-j}^2
\nonumber
\\
&&{} + \sum_{i=1}^{n- 2 \theta\ell-2} \sum
_{j=i+ \theta\ell+2}^{n- \theta\ell} \Sigma_{i,j}^2 +
\sum_{i= \theta\ell+3}^{n- \theta\ell} \sum
_{j=1}^{i- \theta\ell-2} \Sigma_{i,j}^2,
\\
&=& \cases{ O \bigl( n^{-1} \bigr), &\quad if $\nu< (4 \ell-1)/4$,
\vspace*{3pt}
\cr
O \bigl\{ n^{-1} \log(n) \bigr\}, &\quad if $\nu= (4
\ell-1)/4$, \vspace*{3pt}
\cr
O \bigl( n^{-4 \ell+ 4 \nu} \bigr), &\quad if $\nu> (4
\ell-1)/4$,}
\nonumber
\end{eqnarray}
as $n\rightarrow\infty$.
It follows from (\ref{eq:2.5}), (\ref{eq:2.2}), (\ref{eq:2.8}) and the
Borel--Cantelli lemma that $V_{\theta, \ell}/E (V_{\theta, \ell} )
\rightarrow1$ as $n\rightarrow\infty$ almost surely.
Finally,
we observe that for $i, j =1,\ldots, n - \theta\ell$,
\begin{eqnarray*}
&& E \bigl\{ ( \nabla_{\theta, \ell} X_i )^2 (
\nabla_{\theta, \ell} X_j )^2 \bigr\} \label{eq:2.11}
\\
&&\qquad = E \bigl\{ ( \nabla_{\theta, \ell} X_i )^2
\bigr\} E \bigl\{ ( \nabla_{\theta, \ell} X_j )^2 \bigr\}
\\
&&\quad\qquad{} + 2 \Biggl\{ \sum_{k_1 =0}^\ell
\sum_{k_2 =0}^\ell a_{\theta, \ell; i, k_1 }
a_{\theta, \ell; j, k_2 } K( t_{i+ \theta k_1}, t_{j+ \theta k_2}) \Biggr
\}^2,
\nonumber
\\
&& \sum_{i=1}^{n- \theta\ell} \sum
_{j =1}^{n- \theta\ell} \Biggl\{ \sum
_{k_1 =0}^\ell\sum_{k_2 =0}^\ell
a_{\theta, \ell; i, k_1 } a_{\theta, \ell; j, k_2 } K( t_{i+ \theta
k_1}, t_{j+ \theta k_2})
\Biggr\}^2
\nonumber
\\
&&\qquad = \sum_{1\leq i < j \leq n- \theta\ell, j-i \geq\theta\ell+ 2} \Biggl\{ \sum
_{k_1 =0}^\ell\sum_{k_2 =0}^\ell
a_{\theta, \ell; i, k_1 } a_{\theta, \ell; j, k_2 } K( t_{i+ \theta
k_1}, t_{j+ \theta k_2})
\Biggr\}^2
\nonumber
\\
&&\quad\qquad{} + \sum_{1\leq j < i \leq n- \theta\ell, i-j \geq\theta\ell+ 2} \Biggl\{ \sum
_{k_1 =0}^\ell\sum_{k_2 =0}^\ell
a_{\theta, \ell; i, k_1 } a_{\theta, \ell; j, k_2 } K( t_{i+ \theta
k_1}, t_{j+ \theta k_2})
\Biggr\}^2
\nonumber
\\
&&\quad\qquad{} + \sum_{1\leq i, j \leq n- \theta\ell, \llvert i-j\rrvert \leq
\theta\ell+1} \Biggl\{ \sum
_{k_1 =0}^\ell\sum_{k_2 =0}^\ell
a_{\theta, \ell; i, k_1 } a_{\theta, \ell; j, k_2 } K( t_{i+ \theta
k_1}, t_{j+ \theta k_2})
\Biggr\}^2,
\nonumber
\\
&& \sum_{1\leq i, j \leq n- \theta\ell, \llvert i-j\rrvert \leq
\theta\ell+1} \Biggl\{ \sum
_{k_1 =0}^\ell\sum_{k_2 =0}^\ell
a_{\theta, \ell; i, k_1 } a_{\theta, \ell; j, k_2 } K( t_{i+ \theta
k_1}, t_{j+ \theta k_2})
\Biggr\}^2
\nonumber
\\
&&\qquad = O \bigl( n^{4 \ell- 4\nu+1} \bigr),
\nonumber
\\
&& \sum_{1\leq i < j \leq n- \theta\ell, j-i \geq\theta\ell+ 2} \Biggl\{ \sum
_{k_1 =0}^\ell\sum_{k_2 =0}^\ell
a_{\theta, \ell; i, k_1 } a_{\theta, \ell; j, k_2 } K( t_{j+\theta
k_2}, t_{i+ \theta k_1})
\Biggr\}^2
\nonumber
\\
&&\qquad = \sum_{1\leq i < j \leq n- \theta\ell, j-i \geq\theta\ell+ 2} \Biggl\{ \sum
_{k_1 =0}^\ell\sum_{k_2 =0}^\ell
a_{\theta, \ell; i, k_1 } a_{\theta, \ell; j, k_2 }
\nonumber
\\
&&\quad\qquad{} \times\frac{1}{(2 \ell-1)!} \int_{t_j - t_i}^{t_{j +\theta k_2} -
t_{i + \theta k_1}}
(t_{j + \theta k_2} - t_{i+ \theta k_1} -t)^{2 \ell-1} \tilde{K}^{(2
\ell)}
(t) \,d t \Biggr\}^2
\nonumber
\\
&&\qquad = O(1) \sum_{1\leq i < j \leq n- \theta\ell, j-i \geq\theta\ell+
2} \biggl(
\frac{j-i- \theta\ell}{n} \biggr)^{4 \nu -4\ell}
\nonumber
\\
&&\qquad \leq O \bigl(n^2 \bigr) \int_{1/n}^1
s^{4 \nu- 4 \ell} \,d s
\nonumber
\\
&&\qquad = \cases{ O \bigl( n^{4 \ell- 4 \nu+1} \bigr), &\quad if $\nu< (4
\ell-1)/4$, \vspace*{3pt}
\cr
O \bigl\{ n^2 \log(n) \bigr\}, &\quad if
$\nu= (4 \ell-1)/4$, \vspace*{3pt}
\cr
O \bigl( n^2 \bigr), &\quad if
$\nu> (4 \ell-1)/4$,}
\end{eqnarray*}
as $n\rightarrow\infty$. Thus, we conclude that
\[
\frac{ E (V_{\theta, \ell}^2) - \{ E (V_{\theta, \ell}) \}^2}{ \{ E
(V_{\theta, \ell}) \}^2 } = \cases{ O \bigl( n^{-1} \bigr), &\quad if $\nu<
(4 \ell-1)/4$, \vspace*{3pt}
\cr
O \bigl\{ n^{-1} \log(n) \bigr\}, &
\quad if $\nu= (4 \ell-1)/4$, \vspace*{3pt}
\cr
O \bigl( n^{-4 \ell+ 4 \nu} \bigr), &
\quad if $\nu> (4 \ell -1)/4$,}
\]
as $n\rightarrow\infty$.

\item[(b)] Suppose $\nu= \ell$. It follows from Lemma~\ref{la:a.01} that
\begin{eqnarray*}
&& \sum_{i=1}^{n- \theta\ell} E \bigl\{ (
\nabla_{\theta, \ell} X_i)^2 \bigr\}
\nonumber
\\
&&\qquad = \beta_\ell^* \sum_{i=1}^{n- \theta\ell}
\sum_{ k_1 =0}^\ell\sum
_{k_2=0}^\ell a_{\theta, \ell; i, k_1} a_{\theta, \ell; i, k_2}
\nonumber
\\
&&\quad\qquad{} \times(t_{i+ \theta k_2} - t_{i+ \theta k_1} )^{2 \ell}
\log \bigl( \llvert t_{i+ \theta k_2} - t_{i+ \theta k_1}\rrvert \bigr) \bigl\{1 +
o(1) \bigr\}
\nonumber
\\
&&\qquad = - \beta_\ell^* \sum_{i=1}^{n- \theta\ell}
\sum_{ k_1 =0}^\ell\sum
_{k_2=0}^\ell a_{\theta, \ell; i, k_1} a_{\theta, \ell; i, k_2}
(t_{i+ \theta k_2} - t_{i+ \theta k_1} )^{2 \ell} \log(n) \bigl\{1 + o(1)
\bigr\}
\nonumber
\\
&&\qquad \sim (-1)^{\ell+1} \beta_\ell^* (2 \ell)! n \log(n).
\end{eqnarray*}
Next, we observe as in (a) that
\begin{eqnarray*}
&& \sum_{i=1}^{n- \theta\ell} \sum
_{j=1}^{n-\theta\ell} E \bigl\{ ( \nabla_{\theta, \ell}
X_i )^2 ( \nabla_{\theta, \ell} X_j
)^2 \bigr\}
\nonumber
\\
&&\qquad = \sum_{i=1}^{n- \theta\ell} E \bigl\{ (
\nabla_{\theta, \ell} X_i )^2 \bigr\} \sum
_{j=1}^{n- \theta\ell} E \bigl\{ ( \nabla_{\theta, \ell}
X_j )^2 \bigr\}
\nonumber
\\
&&\quad\qquad{} + 2 \sum_{i=1}^{n- \theta\ell} \sum
_{j=1}^{n- \theta\ell} \Biggl\{ \sum
_{k_1 =0}^\ell\sum_{k_2 =0}^\ell
a_{\theta, \ell; i, k_1 } a_{\theta, \ell; j, k_2 } \tilde{K}( t_{i+
\theta k_1} -
t_{j+\theta k_2}) \Biggr\}^2
\nonumber
\\
&&\qquad = \sum_{i=1}^{n- \theta\ell} E \bigl\{ (
\nabla_{\theta, \ell} X_i )^2 \bigr\} \sum
_{j=1}^{n- \theta\ell} E \bigl\{ ( \nabla_{\theta, \ell}
X_j )^2 \bigr\} + O \bigl(n^2 \bigr),
\end{eqnarray*}
as $n\rightarrow\infty$.
Consequently, we conclude that
$\operatorname{Var} \{ V_{\theta, \ell}/ E (V_{\theta, \ell}) \} = O \{ \log
^{-2} (n)\}$ as $n\rightarrow\infty$.

\item[(c)] Suppose $\nu> \ell$. Again using Lemma~\ref{la:a.01}, we have
\begin{eqnarray*}
\sum_{i=1}^{n- \theta\ell} E \bigl\{ (
\nabla_{\theta, \ell} X_i)^2 \bigr\} &\sim& 2
\beta_\ell\sum_{i=1}^{n- \theta\ell} \sum
_{0\leq k_1 <k_2 \leq\ell} a_{\theta, \ell; i, k_1} a_{\theta, \ell;
i, k_2}
(t_{i+ \theta k_1} - t_{i+ \theta k_2} )^{2 \ell}
\nonumber
\\
&=& (-1)^\ell\beta_\ell(2 \ell)! (n - \theta\ell).
\end{eqnarray*}
We observe from Lemmas~\ref{la:a.01} and~\ref{la:a.2} that
\begin{eqnarray*}
\hspace*{-2pt}&& \sum_{i=1}^{n- \theta\ell} \sum
_{j =1}^{n- \theta\ell} \Biggl\{ \sum
_{k_1 =0}^\ell\sum_{k_2=0}^\ell
a_{\theta, \ell; i, k_1 } a_{\theta, \ell; j, k_2 } \tilde{K}( t_{i+
\theta k_1} -
t_{j+ \theta k_2}) \Biggr\}^2
\nonumber
\\
\hspace*{-2pt}&&\qquad = \sum_{i=1}^{n- \theta\ell}
\sum_{j =1}^{n- \theta\ell} \Biggl[ \frac{ \tilde{K}^{(2 \ell)} ( t_i-t_j)
}{(2 \ell-1)!}
\\
\hspace*{-2pt}&&\quad\qquad{}\times \sum_{k_1 =0}^\ell
\sum_{k_2=0}^\ell a_{\theta, \ell; i, k_1 }
a_{\theta, \ell; j, k_2 } \int_0^{t_{i+
\theta k_1} - t_{j+ \theta k_2} - t_i + t_j} t^{2 \ell-1}
\,d t
\nonumber
\\
\hspace*{-2pt}&&\quad\qquad{} + \sum_{k_1 =0}^\ell
\sum_{k_2=0}^\ell\frac{ a_{\theta, \ell; i, k_1 } a_{\theta, \ell; j, k_2
} }{(2 \ell-1)!}
\nonumber
\\
\hspace*{-2pt}&&\quad\qquad{}\times\int_{t_i-t_j}^{t_{i+ \theta k_1} - t_{j+ \theta k_2}}
( t_{i+
\theta k_1} - t_{j+ \theta k_2} - t)^{2 \ell-1} \bigl\{
\tilde{K}^{(2 \ell)} ( t) - \tilde{K}^{(2 \ell)} (t_i -
t_j) \bigr\} \,d t \Biggr]^2
\nonumber
\\
\hspace*{-2pt}&&\qquad= \sum_{i=1}^{n- \theta\ell}
\sum_{j =1}^{n- \theta\ell} \Biggl[ \frac{ \tilde{K}^{(2 \ell)} ( t_i-t_j)
}{(2 \ell)!}
\nonumber
\\
\hspace*{-2pt}&&\quad\qquad{} \times \Biggl\{ \sum_{k_1 =0}^\ell
\sum_{k_2=0}^\ell a_{\theta, \ell; i, k_1 }
a_{\theta, \ell; j, k_2 } ( t_{i+\theta k_1} - t_{j+ \theta k_2} - t_i +
t_j )^{2 \ell} \Biggr\} + o(1) \Biggr]^2
\nonumber
\\
\hspace*{-2pt}&&\qquad = \sum_{i=1}^{n- \theta\ell}
\sum_{j =1}^{n- \theta\ell} \bigl[ \bigl\{
(-1)^\ell \tilde{K}^{(2 \ell)} ( t_i-t_j)
\bigr\}^2 + o(1) \bigr] > c n^2,
\end{eqnarray*}
for some constant $c >0$ since $\tilde{K}^{(2\ell)}(\cdot)$ is a continuous
and not identically 0 function on $(0,1]$. Consequently, $\liminf_{n\rightarrow\infty} \operatorname{Var}\{ V_{\theta, \ell}/E(V_{\theta, \ell})
\}$ equals
\[
\liminf_{n\rightarrow\infty} \frac{ 2}{\{ E (V_{\theta, \ell} ) \}^2} \sum
_{i=1}^{n- \theta\ell} \sum_{j =1}^{n-\theta\ell}
\Biggl\{ \sum_{k_1 =0}^\ell\sum
_{k_2=0}^\ell a_{\theta, \ell; i, k_1 } a_{\theta, \ell; j, k_2 }
\tilde{K}( t_{i+ \theta k_1} - t_{j+\theta k_2}) \Biggr\}^2 >0.
\]
This proves Theorem~\ref{tm:2.1}.\quad\qed
\end{longlist}\noqed
\end{pf*}

%
\begin{la} \label{la:a.4}
Let $p>0$ be an integer and $f_{\theta, \ell} (\cdot), \tilde{f}_{\theta,
2} (\cdot), \bar{f}_{\theta, \ell} (\cdot)$ be as in (\ref{eq:2.101}), (\ref
{eq:2.100}), (\ref{eq:4.00}), respectively. Then
%
%
\begin{eqnarray}
\label{eq:a.78} \lim_{\nu^*\rightarrow p} \frac{ f_{2, \ell} (\nu^*) }{ f_{1, \ell} (\nu
^*) } &=&
\frac{ f_{2, \ell} (p)}{f_{1, \ell} ( p)} \qquad \forall\ell> p,
\nonumber
\\
\lim_{\nu^*\rightarrow1} \frac{ \tilde{f}_{2, 2} (\nu^*) }{ \tilde
{f}_{1, 2} (\nu^*) } &=& \frac{ \tilde{f}_{2, 2} (1)}{ \tilde{f}_{1, 2}
( 1)},
\\
\lim_{\nu^*\rightarrow1} \frac{ \bar{f}_{2, \ell} (\nu^*) }{ \bar
{f}_{1,\ell} (\nu^*) } &=& \frac{ \bar{f}_{2,\ell} (1)}{ \bar{f}_{1,\ell
} ( 1)} \qquad
\forall\ell\in\{1,2\}.
\nonumber
\end{eqnarray}
\end{la}

\begin{pf}
Writing $\nu^* = p + \delta$, we observe from Lemma~\ref
{la:a.01} that for $\theta\in\{1,2\}$,
\begin{eqnarray*}
&& \lim_{\delta\rightarrow0} \frac{1}{2 \delta} \sum
_{i=1}^{n- \theta\ell} \sum_{0\leq k_1 < k_2 \leq l}
a_{\theta, \ell; i, k_1} a_{\theta, \ell; i, k_2} \llvert t_{i+\theta
k_2} - t_{i+\theta k_1}
\rrvert^{2 \nu^*}
\nonumber
\\
&&\qquad = \lim_{\delta\rightarrow0} \frac{1}{2 \delta} \sum
_{i=1}^{n- \theta\ell} \sum_{0\leq k_1 < k_2 \leq l}
a_{\theta, \ell; i, k_1} a_{\theta, \ell; i, k_2} ( t_{i+ \theta k_2} - t_{i+ \theta k_1}
)^{2 p}
\\
&&\quad\qquad{}\times \bigl\{ e^{2 \delta\log( t_{i+ \theta k_2} - t_{i+ \theta k_1}
) } -1 \bigr\}
\nonumber
\\
&&\qquad = \sum_{i=1}^{n- \theta\ell} \sum
_{0\leq k_1 < k_2 \leq l} a_{\theta, \ell; i, k_1} a_{\theta, \ell; i,
k_2} ( t_{i+ \theta k_2}
- t_{i+ \theta k_1} )^{2 p} \log( t_{i+ \theta k_2} - t_{i+ \theta k_1}
)
\end{eqnarray*}
and the first statement in (\ref{eq:a.78}) follows from the definition
of $f_{\theta, \ell} (\cdot)$. The proof of the second statement in (\ref
{eq:a.78}) is similar.
Next, writing $\nu^*= 1 +\delta$ with $\delta\neq0$, we have
\begin{eqnarray*}
&& \lim_{\delta\rightarrow0} \frac{1}{2 \delta} \sum
_{1\leq i_1, i_2\leq n-\theta} \sum_{0\leq k_1, l_1, k_2, l_2 \leq1:
(k_1,k_2)\neq(l_1, l_2)}
c_{\theta, \ell; i_1, i_2}^{k_1, k_2} c_{\theta, \ell; i_1, i_2}^{l_1, l_2}
\nonumber
\\
&&\quad \times \bigl\llVert{\mathbf x}^{i_1 + \theta k_1, i_2 + \theta k_2} - {\mathbf
x}^{i_1 + \theta l_1, i_2 + \theta l_2} \bigr\rrVert^{2 + 2 \delta}
\nonumber
\\
&&\qquad = \lim_{\delta\rightarrow0} \frac{1}{2 \delta} \sum
_{1\leq i_1, i_2\leq n-\theta} \sum_{0\leq k_1, l_1, k_2, l_2 \leq1:
(k_1,k_2)\neq(l_1, l_2)}
c_{\theta, \ell; i_1, i_2}^{k_1, k_2} c_{\theta, \ell; i_1, i_2}^{l_1, l_2}
\nonumber
\\
&&\quad\qquad{}\times \bigl\llVert{\mathbf x}^{i_1 + \theta k_1, i_2 + \theta k_2} - {\mathbf
x}^{i_1 + \theta l_1, i_2 + \theta l_2} \bigr\rrVert^2
\\
&&\quad\qquad{}\times \bigl\{ e^{2 \delta\log(\llVert {\mathbf x}^{i_1 + \theta
k_1, i_2 + \theta k_2} - {\mathbf x}^{i_1 + \theta l_1, i_2 + \theta l_2}
\rrVert )} -1 \bigr\}
\nonumber
\\
&&\qquad = \sum_{1\leq i_1, i_2\leq n-\theta} \sum
_{0\leq k_1, l_1, k_2, l_2
\leq1: (k_1,k_2)\neq(l_1, l_2)} c_{\theta, \ell; i_1, i_2}^{k_1, k_2} c_{\theta, \ell; i_1, i_2}^{l_1, l_2}
\nonumber
\\
&&\quad\qquad{} \times \bigl\llVert{\mathbf x}^{i_1 + \theta k_1, i_2 + \theta k_2} - {\mathbf
x}^{i_1 + \theta l_1, i_2 + \theta l_2} \bigr\rrVert^2
\\
&&\quad\qquad{}\times \log \bigl( n \bigl\llVert{\mathbf x}^{i_1 + \theta k_1, i_2 +
\theta k_2} - {
\mathbf x}^{i_1 + \theta l_1, i_2 + \theta l_2} \bigr\rrVert \bigr).
\end{eqnarray*}
The third statement in (\ref{eq:a.78}) now follows from the definition
of $\bar{f}_{\theta, \ell}(\cdot)$.
\end{pf}

\begin{pf*}{Proof of Theorem~\ref{tm:2.2}} We observe from Lemma~\ref{la:a.1}
that $F_{\ell, n} (\nu^*) \rightarrow2^{2 \nu^* - 2 \ell}$ as
$n\rightarrow\infty$ uniformly over $\nu^* \in[0, M]$.
Suppose with probability 1, $\hat{\nu}_{a, \ell} \rightarrow\nu_a \neq
\nu$ along a subsequence $n_i$ of $n$.
Then
\[
\biggl\{ \frac{ V_{1, \ell} F_{\ell, n_i} (\hat{\nu}_{a, \ell}) }{
V_{2, \ell} } -1 \biggr\}^2 = \biggl\{ \bigl( 1 + o(1)
\bigr) \frac{ F_{\ell, n_i} (\hat{\nu}_{a, \ell}) }{ F_{\ell, n_i} (\nu
) } - 1 \biggr\}^2 \rightarrow \bigl\{
2^{2 (\nu_a - \nu)} - 1 \bigr\}^2 > 0,
\]
as $n_i \rightarrow\infty$ almost surely. This implies that for
sufficiently large $n_i$,
\[
\biggl\{ \frac{ V_{1, \ell} F_{\ell, n_i} (\hat{\nu}_{a, \ell}) }{
V_{2, \ell} } -1 \biggr\}^2 > \frac{1}{2} \bigl
\{ 2^{2 (\nu_a - \nu)} - 1 \bigr\}^2 > \biggl\{ \frac{ V_{1, \ell}
F_{\ell, n_i} (\nu) }{ V_{2, \ell} } -1
\biggr\}^2
\]
almost surely. This contradicts the definition of $\hat{\nu}_{a, \ell
}$.
\end{pf*}

%
\begin{la} \label{la:3.1}
Let $M\geq1$ be an arbitrary but fixed integer. Suppose Conditions~\ref{co2}
and \ref{co4} hold.
Then
\[
d_{i+k_1, i+k_2} = d_{i, i+k_2} - d_{i, i+k_1} + O \bigl(
n^{-3} \bigr),
\]
as $n\rightarrow\infty$ uniformly over $1\leq i\leq n - k_2$ and
$0\leq k_1 \leq k_2\leq M$.
\end{la}

\begin{pf}
For each $t \in[0, L]$, we consider a set of local coordinates of
${\mathbb R}^2$ in a neighbourhood of the point $\gamma(t)$ where the
origin $(0,0)$ of the local coordinates is the point $\gamma(t)$
and the $x$-axis of the local coordinates corresponds to the tangent to
the curve $\gamma$ at $\gamma(t)$. In particular, the vector $\gamma
^{(1)} (t)$ becomes $(1,0)'$.
Since $\gamma$ is a $C^2$-curve, we observe that there exists a
constant $\delta>0$ (independent of $t$), such that under the local
coordinates, $\gamma$ in
a neighbourhood of $\gamma(t)$ can be represented as
$y_t (x) = y_t (0) + y_t^{(1)} (0) x + O( x^2) = O(x^2)$ $\forall0
\leq x \leq\delta$ uniformly in $t\in[0, L]$.
Consequently taking $t = t_i$ and writing $(x_{i+k}, y_{t_i}
(x_{i+k}))$ for the point $\gamma( t_{i+k})$ in the local coordinates,
we obtain $x_i = 0 = y_{t_i} ( 0) = y_{t_i}^{(1)} ( 0)$ and
\begin{eqnarray*}
d_{i+k_1, i+ k_2} &=& \sqrt{ (x_{i+k_2} - x_{i+k_1})^2
+ \bigl\{ (y_{t_i} ( x_{i+k_2}) - y_{t_i} (
x_{i+k_1} ) \bigr\}^2 }
\nonumber
\\
&=& (x_{i+k_2} - x_{i+k_1} ) \bigl\{1 + O \bigl(
n^{-2} \bigr) \bigr\}
\nonumber
\\
&=& d_{i, i+k_2} - d_{i, i+k_1} + O \bigl(n^{-3} \bigr),
\end{eqnarray*}
as $n\rightarrow\infty$ uniformly over $1\leq i\leq n- k_2$ and $0\leq
k_1 \leq k_2 \leq M$.
\end{pf}

%
\begin{la} \label{la:3.9}
Let $\theta, \ell\in\{1,2 \}$, $\tilde{f}_{\theta, \ell} (\cdot)$ be as
in (\ref{eq:2.100}) and $H_\ell(\cdot)$ be as in Lemma~\ref{la:a.1}.
Suppose Conditions~\ref{co2},~\ref{co3} and~\ref{co4} hold.
Then for $\nu\in(0, \ell)$,
\[
\tilde{f}_{\theta, \ell} (\nu) = \frac{2 \beta^*_\nu n^{2\ell+1 -2\nu
} H_\ell(\nu) }{ \theta^{2 \ell-2\nu} L^{2 \ell+1 -2\nu}} \int_0^L
\bigl\llVert(\gamma\circ\varphi)^{(1)} ( t ) \bigr\rrVert
^{2 \nu- 2 \ell} \,dt + O \bigl(n^{2 \ell-2 \nu} \bigr),
\]
and hence
$\tilde{F}_{\ell, n} (\nu^*) = 2^{2\nu^* - 2\ell} + O(n^{-1})$ as
$n\rightarrow\infty$
uniformly over $0 \leq\nu^* \leq\ell$ where $\tilde{F}_{\ell, n}(\cdot)$
is as in (\ref{eq:9.69}).
\end{la}

\begin{pf}
We observe that for $0\leq k \leq4$,
\begin{eqnarray*}
d_{i, i+k} &=& \biggl[ \biggl\{ \gamma_1 \circ\varphi \biggl(
\frac{ L ( i + k -1) }{ n-1} \biggr) - \gamma_1 \circ\varphi \biggl(
\frac{ L ( i -1) }{ n-1} \biggr) \biggr\}^2
\nonumber
\\
&&{} + \biggl\{ \gamma_2 \circ\varphi \biggl( \frac{ L ( i + k -1) }{ n-1}
\biggr) - \gamma_2 \circ\varphi \biggl( \frac{ L ( i -1) }{ n-1} \biggr)
\biggr\}^2 \biggr]^{1/2}
\nonumber
\\
&=& \frac{ k L }{ n-1}
\biggl\llVert(\gamma \circ\varphi)^{(1)} \biggl( \frac{ L ( i -1) }{ n-1}
\biggr) \biggr\rrVert+ O \bigl(n^{-2} \bigr),
\end{eqnarray*}
as $n\rightarrow\infty$ uniformly over $1\leq i\leq n-k$.
For $\theta\in\{1,2\}$ and $\nu\in(0,1)$,
\begin{eqnarray*}
&& \sum_{i=1}^{n-\theta} b_{\theta, 1; i, 0}
b_{\theta, 1; i, 1} d_{i, i+ \theta}^{2 \nu} 
\nonumber
\\
&&\qquad = - \sum_{i=1}^{n-\theta} \biggl(
\frac{ \theta L}{n-1} \biggr)^{2\nu-2} \biggl\llVert(\gamma\circ
\varphi)^{(1)} \biggl( \frac{ L ( i -1) }{ n-1} \biggr) \biggr\rrVert
^{2\nu-2} + O \bigl(n^{2 - 2 \nu} \bigr)
\nonumber
\\
&&\qquad = - \theta^{2\nu-2} L^{2\nu-3} n^{3 -2 \nu} \int
_0^L \bigl\llVert(\gamma\circ
\varphi)^{(1)} ( t ) \bigr\rrVert^{2\nu-2} \,dt + O \bigl(
n^{2-2 \nu} \bigr),
\end{eqnarray*}
as $n\rightarrow\infty$.
We observe from Lemma~\ref{la:3.1} that for $\theta\in\{1,2\}$ and
$\nu\in(0,2)\setminus\{1\}$,
\begin{eqnarray*}
&& \sum_{0 \leq k_1 < k_2\leq2} b_{\theta, 2; i, k_1} b_{\theta, 2; i, k_2}
d_{i+ \theta k_1, i+ \theta k_2}^{2 \nu}
\nonumber
\\
&&\qquad = - \frac{ 4 d_{i, i+ \theta}^{2 \nu-2} }{ d_{i,i+2 \theta} (d_{i,
i+2 \theta} - d_{i, i+ \theta}) } + \frac{ 4 d_{i, i+ 2 \theta}^{2 \nu
-2} }{ d_{i,i+ \theta} (d_{i, i+2 \theta} - d_{i, i+\theta}) }
\nonumber
\\
&&\quad\qquad{} - \frac{ 4 \{ d_{i, i+2 \theta} - d_{i,i+ \theta} + O(n^{-3}) \}^{2
\nu} }{ d_{i,i+ \theta} d_{i, i+2 \theta}(d_{i, i+2 \theta} - d_{i,
i+\theta})^2 }
\nonumber
\\
&&\qquad = - 2 \biggl(\frac{ \theta L}{n-1} \biggr)^{2\nu-4} \biggl\llVert (
\gamma \circ\varphi)^{(1)} \biggl( \frac{ L ( i -1) }{ n-1} \biggr) \biggr
\rrVert ^{2 \nu-4}
\nonumber
\\
&&\quad\qquad{} + 2^{2 \nu} \biggl(\frac{ \theta L}{n-1}
\biggr)^{2\nu-4} \biggl\llVert(\gamma\circ\varphi)^{(1)} \biggl(
\frac{ L ( i -1) }{ n-1} \biggr) \biggr\rrVert^{2 \nu-4}
\nonumber
\\
&&\quad\qquad{} - 2 \biggl(\frac{ \theta L}{n-1} \biggr)^{2\nu-4} \biggl
\llVert (\gamma \circ\varphi)^{(1)} \biggl( \frac{ L ( i -1) }{ n-1} \biggr)
\biggr\rrVert ^{2 \nu-4} + O \bigl(n^{3-2 \nu} \bigr)
\nonumber
\\
&&\qquad = \bigl(2^{2 \nu} - 4 \bigr) \biggl(\frac{ \theta L}{n-1}
\biggr)^{2\nu-4} \biggl\llVert(\gamma\circ\varphi)^{(1)} \biggl(
\frac{ L ( i -1) }{ n-1} \biggr) \biggr\rrVert^{2 \nu-4} + O \bigl(n^{3-2 \nu}
\bigr),
\end{eqnarray*}
as $n\rightarrow\infty$ uniformly over $1\leq i\leq n-2 \theta$.
Consequently,
\begin{eqnarray*}
&& \sum_{i=1}^{n-2 \theta} \sum
_{0 \leq k_1 < k_2\leq2} b_{\theta, 2; i, k_1} b_{\theta, 2; i, k_2}
d_{i+ \theta k_1, i+ \theta k_2}^{2 \nu}
\nonumber
\\
&&\qquad = \bigl(2^{2 \nu} - 4 \bigr) \biggl(\frac{ \theta L}{n-1}
\biggr)^{2\nu-4} \sum_{i=1}^{n-2 \theta}
\biggl\llVert(\gamma\circ\varphi)^{(1)} \biggl( \frac{ L ( i -1) }{ n-1}
\biggr) \biggr\rrVert^{2 \nu-4} + O \bigl(n^{4-2 \nu} \bigr)
\nonumber
\\
&&\qquad = \bigl(2^{2 \nu} -
4 \bigr) \theta^{2 \nu-4} L^{2 \nu-5} n^{5-2\nu} \int
_0^L \bigl\llVert(\gamma\circ
\varphi)^{(1)} ( t ) \bigr\rrVert ^{2 \nu-4} \,dt + O
\bigl(n^{4-2 \nu} \bigr)
\end{eqnarray*}
as $n\rightarrow\infty$.
We observe from Lemma~\ref{la:3.1} that for $\theta\in\{1,2\}$,
\begin{eqnarray*}
&& \sum_{0 \leq k_1 < k_2\leq2} b_{\theta, 2; i, k_1} b_{\theta, 2; i, k_2}
d_{i+ \theta k_1, i+ \theta k_2}^2 \log( d_{i+ \theta k_1, i+ \theta
k_2} )
\nonumber
\\
&&\qquad = \frac{ 4 d_{i, i+ \theta}^2 \log( d_{i, i+ \theta} ) }{ (- d_{i,
i+\theta}) ( - d_{i, i+2 \theta}) d_{i, i+ \theta} (d_{i, i+ \theta} -
d_{i, i+2 \theta}) }
\nonumber
\\
&&\quad\qquad{} + \frac{ 4 d_{i, i+ 2\theta}^2 \log(d_{i, i+ 2 \theta} ) }{ ( -
d_{i, i+ \theta}) ( - d_{i, i+2 \theta}) d_{i, i+2 \theta} (d_{i, i+2
\theta} - d_{i, i+\theta}) }
\nonumber
\\
&&\quad\qquad{} + \frac{ 4 \{ d_{i, i+2 \theta} - d_{i,i+ \theta} + O(n^{-3}) \}^2
\log\{ \llvert d_{i, i+2 \theta} - d_{i,i+ \theta}\rrvert +
O(n^{-3}) \} }{
d_{i,i+ \theta} (d_{i,i+ \theta} - d_{i, i+2 \theta}) d_{i, i+2 \theta}
(d_{i, i+2 \theta} - d_{i, i+ \theta}) }
\nonumber
\\
&&\qquad = \frac{4}{ d_{i,i+\theta}^2} \biggl\{ \biggl( \frac{ d_{i,i+2 \theta
} }{ d_{i, i+\theta} } - 1
\biggr)^{-1} \log \biggl( \frac{d_{i, i+ 2 \theta} }{ d_{i, i+ \theta}
} \biggr) -
\frac{ d_{i,i+\theta} }{ d_{i, i+2 \theta} } \log \biggl( \biggl\llvert \frac{d_{i, i+2 \theta} }{ d_{i, i+ \theta} } -1 \biggr\rrvert
\biggr) \biggr\}
\nonumber
\\
&&\quad\qquad{} + O \bigl\{ \log(n) \bigr\}
\nonumber
\\
&&\qquad = \frac{4 \log(2) n^2}{\theta^2 L^2} \biggl\llVert(\gamma\circ \varphi)^{(1)}
\biggl( \frac{L( i-1) }{n-1}\biggr) \biggr\rrVert^{-2} + O(n),
\end{eqnarray*}
as $n\rightarrow\infty$ uniformly over $1\leq i\leq n-2 \theta$. Consequently,
\begin{eqnarray*}
&& \sum_{i=1}^{n-2 \theta} \sum
_{0 \leq k_1 < k_2\leq2} b_{\theta, 2; i, k_1} b_{\theta, 2; i, k_2}
d_{i+ \theta k_1, i+ \theta k_2}^2 \log( d_{i+ \theta k_1, i+ \theta
k_2} )
\nonumber
\\
&&\qquad = \frac{4 \log(2) n^2}{\theta^2 L^2} \sum_{i=1}^{n-2 \theta}
\biggl\llVert(\gamma\circ\varphi)^{(1)} \biggl( \frac{L( i-1) }{n-1} \biggr)\biggr
\rrVert^{-2} + O \bigl(n^2 \bigr)
\nonumber
\\
&&\qquad = \frac{4 \log(2) n^3 }{\theta^2 L^3} \int_0^L \bigl
\llVert (\gamma\circ\varphi)^{(1)} ( t) \bigr\rrVert^{-2} \,dt
+ O \bigl(n^2 \bigr),
\end{eqnarray*}
as $n\rightarrow\infty$. This proves the lemma.
\end{pf}

\begin{pf*}{Proof of Lemma~\ref{la:7.1}}
We shall prove the case $\ell=1$. The proof for the case $\ell=2$ is
similar and is omitted.
We observe from (\ref{eq:8.1}) and (\ref{eq:8.2}) that for $\theta\in\{
1, 2\}$,
\begin{eqnarray*}
\alpha_{\theta; i_1, i_2}^{1, 1} \bigl( x_1^{i_1+ \theta, i_2} -
x_1^{i_1, i_2} \bigr) + \alpha_{\theta; i_1, i_2}^{1, 2}
\bigl( x_1^{i_1, i_2+ \theta} - x_1^{i_1, i_2} \bigr) &=&
1,
\nonumber
\\
\alpha_{\theta; i_1, i_2}^{1, 1} \bigl( x_2^{i_1+ \theta, i_2} -
x_2^{i_1, i_2} \bigr) +\alpha_{\theta; i_1,i_2}^{1, 2}
\bigl(x_2^{i_1, i_2+ \theta} - x_2^{i_1, i_2} \bigr) &=&
0,
\nonumber
\\
\beta_{\theta; i_1, i_2}^{1, 1} \bigl( x_1^{i_1+ \theta, i_2} -
x_1^{i_1+ \theta, i_2+ \theta} \bigr) + \beta_{\theta; i_1,i_2}^{1, 2}
\bigl(x_1^{i_1, i_2+ \theta} - x_1^{i_1+ \theta, i_2+ \theta} \bigr) &=&
1,
\nonumber
\\
\beta_{\theta; i_1, i_2}^{1, 1} \bigl( x_2^{i_1+ \theta, i_2} -
x_2^{i_1+ \theta, i_2+ \theta} \bigr) + \beta_{\theta; i_1, i_2}^{1, 2}
\bigl(x_2^{i_1, i_2+ \theta} - x_2^{i_1+ \theta, i_2+ \theta} \bigr) &=&
0.
\nonumber
\\
\end{eqnarray*}
Lemma~\ref{la:7.1} is a consequence of (\ref{eq:8.3}) and the above
four equalities.
\end{pf*}

%
\begin{la} \label{la:7.11}
For $0 < \nu< 2$ and $\ell\in\{1,2\}$, define ${\mathcal J}_{\nu,
\ell}: [0,1]^2 \rightarrow{\mathbb R}$ by
\begin{eqnarray*}
{\mathcal J}_{\nu, \ell} (u,v) &=& \frac{ \{ \varphi_{\ell^c}^{(0,1)} (
u, v) - \varphi_{\ell^c}^{(1,0)} ( u, v) \}^2 }{
\{ \varphi_1^{(1,0)} ( u, v) \varphi_2^{(0,1)} ( u, v)
- \varphi_1^{(0,1)} ( u, v) \varphi_2^{(1,0)} ( u, v) \}^2 }
\nonumber
\\
&&{} \times \Biggl\{ - 2 G_\nu \Biggl( \Biggl[ \sum
_{j=1}^2 \varphi_j^{(0,1)} (
u,v)^2 \Biggr]^{1/2} \Biggr) - 2 G_\nu \Biggl(
\Biggl[ \sum_{j=1}^2 \varphi_j^{(1,0)}
( u, v)^2 \Biggr]^{1/2} \Biggr)
\nonumber
\\
&&{} + G_\nu \Biggl( \Biggl[ \sum_{j=1}^2
\bigl\{ \varphi_j^{(1,0)} ( u, v) + \varphi_j^{(0,1)}
( u, v) \bigr\}^2 \Biggr]^{1/2} \Biggr)
\nonumber
\\
&&{} + G_\nu \Biggl( \Biggl[ \sum_{j=1}^2
\bigl\{ \varphi_j^{(0,1)} ( u, v) - \varphi_j^{(1,0)}
( u, v) \bigr\}^2 \Biggr]^{1/2} \Biggr) \Biggr\},
\end{eqnarray*}
where $\ell^c = 2$ if $\ell=1$ and $\ell^c = 1$ if $\ell=2$. Then $\min_{0\leq u, v\leq1} {\mathcal J}_{\nu, \ell} (u,v) \geq0$ if $\nu\in
[1, 2)$,
and $\max_{0\leq u, v\leq1} {\mathcal J}_{\nu, \ell} (u,v) \leq0$ if
$\nu\in(0, 1)$.
\end{la}

%
\begin{la} \label{la:7.13}
Let $\theta, \ell\in\{1,2 \}$ and $\bar{f}_{\theta, \ell} (\cdot), \bar
{F}_{n,\ell}(\cdot)$, be as in (\ref{eq:4.00}), (\ref{eq:10.1}), respectively.
Suppose Condition~\ref{co5} holds and that ${\mathcal J}_{\nu, \ell} (\cdot,\cdot)$, as
in Lemma~\ref{la:7.11}, is not identically $0$ on $[0,1]^2$.
Then
\[
\bar{f}_{\theta, \ell} (\nu) = 2 \beta_\nu^* \theta^{2\nu-2}
n^{4-2\nu} \int_0^1 \int
_0^1 {\mathcal J}_{\nu, \ell} (u, v) \,du
\,dv + O \bigl(n^{3 - 2 \nu} \bigr),
\]
and hence $\bar{F}_{n,\ell} (\nu^*) = 2^{2 \nu^* -2} + O(n^{-1})$ as $n
\rightarrow\infty$ uniformly over $0 \leq\nu^* \leq2$.
\end{la}

We refer the reader to \cite{l2} for the proofs of Lemmas~\ref{la:7.11}
and~\ref{la:7.13}.
\end{appendix}

\section*{Acknowledgements}
I would like to thank Professor Runze Li, an Associate Editor and four
referees for their comments and suggestions
that led to numerous improvements in the initial results of this article.

\begin{supplement}[id=suppA]
\stitle{Proofs and other technical details}
\slink[doi]{10.1214/15-AOS1365SUPP} 
\sdatatype{.pdf}
\sfilename{aos1365\_supp.pdf}
\sdescription{The supplemental article \cite{l2} contains the proofs of
Theorems~\ref{tm:4.1},~\ref{tm:7.1},~\ref{tm:10.1}, Lemmas~\ref
{la:7.11},~\ref{la:7.13} and Corollary~\ref{cy:7.1}.}
\end{supplement}

%

\printaddresses

\begin{thebibliography}{28}

\bibitem{ap1}
%
\begin{barticle}[mr]
\bauthor{\bsnm{Adler},~\bfnm{Robert~J.}\binits{R.~J.}} \AND
\bauthor{\bsnm{Pyke},~\bfnm{Ron}\binits{R.}}
(\byear{1993}).
\btitle{Uniform quadratic variation for {G}aussian processes}.
\bjournal{Stochastic Process. Appl.}
\bvolume{48}
\bpages{191--209}.
\bid{doi={10.1016/0304-4149(93)90044-5}, issn={0304-4149}, mr={1244542}}
\end{barticle}
%
\bptok{imsref}%
\endbibitem

\bibitem{ac1}
%
\begin{barticle}[mr]
\bauthor{\bsnm{Anderes},~\bfnm{Ethan}\binits{E.}} \AND
\bauthor{\bsnm{Chatterjee},~\bfnm{Sourav}\binits{S.}}
(\byear{2009}).
\btitle{Consistent estimates of deformed isotropic {G}aussian random
fields on the plane}.
\bjournal{Ann. Statist.}
\bvolume{37}
\bpages{2324--2350}.
\bid{doi={10.1214/08-AOS647}, issn={0090-5364}, mr={2543694}}
\end{barticle}
%
\bptok{imsref}%
\endbibitem

\bibitem{as1}
%
\begin{barticle}[mr]
\bauthor{\bsnm{Anderes},~\bfnm{Ethan~B.}\binits{E.~B.}} \AND
\bauthor{\bsnm{Stein},~\bfnm{Michael~L.}\binits{M.~L.}}
(\byear{2008}).
\btitle{Estimating deformations of isotropic {G}aussian random fields
on the plane}.
\bjournal{Ann. Statist.}
\bvolume{36}
\bpages{719--741}.
\bid{doi={10.1214/009053607000000893}, issn={0090-5364}, mr={2396813}}
\end{barticle}
%
\bptok{imsref}%
\endbibitem

\bibitem{b1}
%
\begin{barticle}[mr]
\bauthor{\bsnm{Begyn},~\bfnm{Arnaud}\binits{A.}}
(\byear{2005}).
\btitle{Quadratic variations along irregular subdivisions for
{G}aussian processes}.
\bjournal{Electron. J. Probab.}
\bvolume{10}
\bpages{691--717 (electronic)}.
\bid{doi={10.1214/EJP.v10-245}, issn={1083-6489}, mr={2164027}}
\bptnote{check pages}%
\end{barticle}
%
\bptok{imsref}%
\endbibitem

\bibitem{bcij1}
%
\begin{barticle}[mr]
\bauthor{\bsnm{Benassi},~\bfnm{Albert}\binits{A.}},
\bauthor{\bsnm{Cohen},~\bfnm{Serge}\binits{S.}},
\bauthor{\bsnm{Istas},~\bfnm{Jacques}\binits{J.}} \AND
\bauthor{\bsnm{Jaffard},~\bfnm{St{\'e}phane}\binits{S.}}
(\byear{1998}).
\btitle{Identification of filtered white noises}.
\bjournal{Stochastic Process. Appl.}
\bvolume{75}
\bpages{31--49}.
\bid{doi={10.1016/S0304-4149(97)00123-3}, issn={0304-4149}, mr={1629014}}
\end{barticle}
%
\bptok{imsref}%
\endbibitem

\bibitem{cw1}
%
\begin{barticle}[mr]
\bauthor{\bsnm{Chan},~\bfnm{Grace}\binits{G.}} \AND
\bauthor{\bsnm{Wood},~\bfnm{Andrew~T.~A.}\binits{A.~T.~A.}}
(\byear{2000}).
\btitle{Increment-based estimators of fractal dimension for
two-dimensional surface data}.
\bjournal{Statist. Sinica}
\bvolume{10}
\bpages{343--376}.
\bid{issn={1017-0405}, mr={1769748}}
\end{barticle}\vadjust{\goodbreak}
%
\bptok{imsref}%
\endbibitem

\bibitem{cd1}
%
\begin{bbook}[mr]
\bauthor{\bsnm{Chil{\`e}s},~\bfnm{Jean-Paul}\binits{J.-P.}} \AND
\bauthor{\bsnm{Delfiner},~\bfnm{Pierre}\binits{P.}}
(\byear{1999}).
\btitle{Geostatistics. Modeling Spatial Uncertainty}.
\bpublisher{Wiley},
\blocation{New York}.
\bid{doi={10.1002/9780470316993}, mr={1679557}}
\end{bbook}
%
\bptok{imsref}%
\endbibitem

\bibitem{cgpp1}
%
\begin{barticle}[mr]
\bauthor{\bsnm{Cohen},~\bfnm{Serge}\binits{S.}},
\bauthor{\bsnm{Guyon},~\bfnm{Xavier}\binits{X.}},
\bauthor{\bsnm{Perrin},~\bfnm{Olivier}\binits{O.}} \AND
\bauthor{\bsnm{Pontier},~\bfnm{Monique}\binits{M.}}
(\byear{2006}).
\btitle{Singularity functions for fractional processes: Application to
the fractional {B}rownian sheet}.
\bjournal{Ann. Inst. Henri Poincar\'e Probab. Stat.}
\bvolume{42}
\bpages{187--205}.
\bid{doi={10.1016/j.anihpb.2005.03.002}, issn={0246-0203}, mr={2199797}}
\end{barticle}
%
\bptok{imsref}%
\endbibitem

\bibitem{ch1}
%
\begin{barticle}[mr]
\bauthor{\bsnm{Constantine},~\bfnm{A.~G.}\binits{A.~G.}} \AND
\bauthor{\bsnm{Hall},~\bfnm{Peter}\binits{P.}}
(\byear{1994}).
\btitle{Characterizing surface smoothness via estimation of effective
fractal dimension}.
\bjournal{J. Roy. Statist. Soc. Ser. B}
\bvolume{56}
\bpages{97--113}.
\bid{issn={0035-9246}, mr={1257799}}
\end{barticle}
%
\bptok{imsref}%
\endbibitem

\bibitem{c1}
%
\begin{bbook}[mr]
\bauthor{\bsnm{Cressie},~\bfnm{Noel~A.~C.}\binits{N.~A.~C.}}
(\byear{1991}).
\btitle{Statistics for Spatial Data}.
\bpublisher{Wiley},
\blocation{New York}.
\bid{mr={1127423}}
\end{bbook}
%
\bptok{imsref}%
\endbibitem

\bibitem{gm1}
%
\begin{bbook}[mr]
\bauthor{\bsnm{Grenander},~\bfnm{Ulf}\binits{U.}} \AND
\bauthor{\bsnm{Miller},~\bfnm{Michael~I.}\binits{M.~I.}}
(\byear{2007}).
\btitle{Pattern Theory: From Representation to Inference}.
\bpublisher{Oxford Univ. Press},
\blocation{Oxford}.
\bid{mr={2285439}}
\end{bbook}
%
\bptok{imsref}%
\endbibitem

\bibitem{gl1}
%
\begin{barticle}[mr]
\bauthor{\bsnm{Guyon},~\bfnm{Xavier}\binits{X.}} \AND
\bauthor{\bsnm{Le{\'o}n},~\bfnm{Jos{\'e}}\binits{J.}}
(\byear{1989}).
\btitle{Convergence en loi des {$H$}-variations d'un processus gaussien
stationnaire sur {${\mathbf R}$}}.
\bjournal{Ann. Inst. Henri Poincar\'e Probab. Stat.}
\bvolume{25}
\bpages{265--282}.
\bid{issn={0246-0203}, mr={1023952}}
\end{barticle}
%
\bptok{imsref}%
\endbibitem

\bibitem{hw1}
%
\begin{barticle}[mr]
\bauthor{\bsnm{Hall},~\bfnm{Peter}\binits{P.}} \AND
\bauthor{\bsnm{Wood},~\bfnm{Andrew}\binits{A.}}
(\byear{1993}).
\btitle{On the performance of box-counting estimators of fractal dimension}.
\bjournal{Biometrika}
\bvolume{80}
\bpages{246--252}.
\bid{doi={10.1093/biomet/80.1.246}, issn={0006-3444}, mr={1225230}}
\end{barticle}
%
\bptok{imsref}%
\endbibitem

\bibitem{hw2}
%
\begin{barticle}[mr]
\bauthor{\bsnm{Hanson},~\bfnm{D.~L.}\binits{D.~L.}} \AND
\bauthor{\bsnm{Wright},~\bfnm{F.~T.}\binits{F.~T.}}
(\byear{1971}).
\btitle{A bound on tail probabilities for quadratic forms in
independent random variables}.
\bjournal{Ann. Math. Statist.}
\bvolume{42}
\bpages{1079--1083}.
\bid{issn={0003-4851}, mr={0279864}}
\end{barticle}
%
\bptok{imsref}%
\endbibitem

\bibitem{hj1}
%
\begin{bbook}[mr]
\bauthor{\bsnm{Horn},~\bfnm{Roger~A.}\binits{R.~A.}} \AND
\bauthor{\bsnm{Johnson},~\bfnm{Charles~R.}\binits{C.~R.}}
(\byear{1985}).
\btitle{Matrix Analysis}.
\bpublisher{Cambridge Univ. Press},
\blocation{Cambridge}.
\bid{doi={10.1017/CBO9780511810817}, mr={0832183}}
\end{bbook}
%
\bptok{imsref}%
\endbibitem

\bibitem{isz1}
%
\begin{barticle}[mr]
\bauthor{\bsnm{Im},~\bfnm{Hae~Kyung}\binits{H.~K.}},
\bauthor{\bsnm{Stein},~\bfnm{Michael~L.}\binits{M.~L.}} \AND
\bauthor{\bsnm{Zhu},~\bfnm{Zhengyuan}\binits{Z.}}
(\byear{2007}).
\btitle{Semiparametric estimation of spectral density with irregular
observations}.
\bjournal{J. Amer. Statist. Assoc.}
\bvolume{102}
\bpages{726--735}.
\bid{doi={10.1198/016214507000000220}, issn={0162-1459}, mr={2381049}}
\end{barticle}
%
\bptok{imsref}%
\endbibitem

\bibitem{il1}
%
\begin{barticle}[mr]
\bauthor{\bsnm{Istas},~\bfnm{Jacques}\binits{J.}} \AND
\bauthor{\bsnm{Lang},~\bfnm{Gabriel}\binits{G.}}
(\byear{1997}).
\btitle{Quadratic variations and estimation of the local H\"older index
of a {G}aussian process}.
\bjournal{Ann. Inst. Henri Poincar\'e Probab. Stat.}
\bvolume{33}
\bpages{407--436}.
\bid{doi={10.1016/S0246-0203(97)80099-4}, issn={0246-0203}, mr={1465796}}
\bptnote{check pages}%
\end{barticle}
%
\bptok{imsref}%
\endbibitem

\bibitem{jm1}
%
\begin{barticle}[mr]
\bauthor{\bsnm{Joshi},~\bfnm{Sarang~C.}\binits{S.~C.}} \AND
\bauthor{\bsnm{Miller},~\bfnm{Michael~I.}\binits{M.~I.}}
(\byear{2000}).
\btitle{Landmark matching via large deformation diffeomorphisms}.
\bjournal{IEEE Trans. Image Process.}
\bvolume{9}
\bpages{1357--1370}.
\bid{doi={10.1109/83.855431}, issn={1057-7149}, mr={1808275}}
\end{barticle}
%
\bptok{imsref}%
\endbibitem

\bibitem{kw1}
%
\begin{barticle}[mr]
\bauthor{\bsnm{Kent},~\bfnm{John~T.}\binits{J.~T.}} \AND
\bauthor{\bsnm{Wood},~\bfnm{Andrew~T.~A.}\binits{A.~T.~A.}}
(\byear{1997}).
\btitle{Estimating the fractal dimension of a locally self-similar
{G}aussian process by using increments}.
\bjournal{J. Roy. Statist. Soc. Ser. B}
\bvolume{59}
\bpages{679--699}.
\bid{issn={0035-9246}, mr={1452033}}
\end{barticle}
%
\bptok{imsref}%
\endbibitem

\bibitem{kg1}
%
\begin{barticle}[mr]
\bauthor{\bsnm{Klein},~\bfnm{Ruben}\binits{R.}} \AND
\bauthor{\bsnm{Gin{\'e}},~\bfnm{Evarist}\binits{E.}}
(\byear{1975}).
\btitle{On quadratic variation of processes with {G}aussian increments}.
\bjournal{Ann. Probab.}
\bvolume{3}
\bpages{716--721}.
\bid{mr={0378070}}
\bptnote{check pages}%
\end{barticle}
%
\bptok{imsref}%
\endbibitem

\bibitem{k1}
%
\begin{bbook}[mr]
\bauthor{\bsnm{Knuth},~\bfnm{Donald~E.}\binits{D.~E.}}
(\byear{1997}).
\btitle{The Art of Computer Programming},
\bedition{3rd} ed.
\bseries{Fundamental Algorithms}
\bvolume{1}.
\bpublisher{Addison-Wesley},
\blocation{Reading, MA}.
\bid{mr={3077152}}
\bptnote{check year}%
\end{bbook}
%
\bptok{imsref}%
\endbibitem

\bibitem{l1}
%
\begin{barticle}[mr]
\bauthor{\bsnm{L{\'e}vy},~\bfnm{Paul}\binits{P.}}
(\byear{1940}).
\btitle{Le mouvement brownien plan}.
\bjournal{Amer. J. Math.}
\bvolume{62}
\bpages{487--550}.
\bid{issn={0002-9327}, mr={0002734}}
\end{barticle}
%
\bptok{imsref}%
\endbibitem

\bibitem{l2}
%
\begin{bmisc}[author]
\bauthor{\bsnm{Loh},~\bfnm{W.-L.}\binits{W.-L.}}
(\byear{2015}).
\bhowpublished{Supplement to ``Estimating the smoothness of a Gaussian
random field from irregularly spaced data via higher-order quadratic
variations''.
DOI:\doiurl{10.1214/15-AOS1365SUPP}}.
\bptok{imsref}%
\end{bmisc}
%
\bptok{imsref}%
\endbibitem

\bibitem{m1}
%
\begin{bbook}[auto:parserefs-M02]
\bauthor{\bsnm{Matheron},~\bfnm{G.}\binits{G.}}
(\byear{1971}).
\btitle{The Theory of Regionalized Variables and Its Applications}.
\bpublisher{Ecole des Mines},
\blocation{Fontainebleau}.
\end{bbook}
%
\bptok{imsref}%
\endbibitem

\bibitem{s1}
%
\begin{bbook}[mr]
\bauthor{\bsnm{Stein},~\bfnm{Michael~L.}\binits{M.~L.}}
(\byear{1999}).
\btitle{Interpolation of Spatial Data: Some Theory for Kriging}.
\bpublisher{Springer},
\blocation{New York}.
\bid{doi={10.1007/978-1-4612-1494-6}, mr={1697409}}
\end{bbook}
%
\bptok{imsref}%
\endbibitem

\bibitem{s2}
%
\begin{barticle}[mr]
\bauthor{\bsnm{Stein},~\bfnm{Michael~L.}\binits{M.~L.}}
(\byear{2012}).
\btitle{Simulation of {G}aussian random fields with one derivative}.
\bjournal{J. Comput. Graph. Statist.}
\bvolume{21}
\bpages{155--173}.
\bid{doi={10.1198/jcgs.2010.10069}, issn={1061-8600}, mr={2913361}}
\end{barticle}
%
\bptok{imsref}%
\endbibitem

\bibitem{s3}
%
\begin{barticle}[auto:parserefs-M02]
\bauthor{\bsnm{Sylvester},~\bfnm{J.~J.}\binits{J.~J.}}
(\byear{1857}).
\btitle{On the partition of numbers}.
\bjournal{Quart. J. Math.}
\bvolume{1}
\bpages{141--152}.
\end{barticle}
%
\bptok{imsref}%
\endbibitem

\bibitem{wc1}
%
\begin{barticle}[mr]
\bauthor{\bsnm{Wood},~\bfnm{Andrew~T.~A.}\binits{A.~T.~A.}} \AND
\bauthor{\bsnm{Chan},~\bfnm{Grace}\binits{G.}}
(\byear{1994}).
\btitle{Simulation of stationary {G}aussian processes in {$[0,1]\sp d$}}.
\bjournal{J. Comput. Graph. Statist.}
\bvolume{3}
\bpages{409--432}.
\bid{doi={10.2307/1390903}, issn={1061-8600}, mr={1323050}}
\end{barticle}
%
\bptok{imsref}%
\endbibitem

\end{thebibliography}
\end{document}